\documentclass[aop,MSNbibl,number,seceqn,dvips]{arximspdf}
\usepackage{graphicx}

\doi{10.1214/11-AOP727} 
\volume{41}
\issue{3B}
\pubyear{2013}
\firstpage{1767}
\lastpage{1805}

\makeatletter
\newcommand{\eqref}[1]{(\ref{#1})}
\newtheorem{ittheorem}{Theorem}[section]
\newtheorem{itlemma}[ittheorem]{Lemma}
\newproclaim{itdefinition}[ittheorem]{Definition}
\newproclaim{itremark}[ittheorem]{Remark}
\newtheorem{itconjecture}[ittheorem]{Conjecture}
\newtheorem{itcorollary}[ittheorem]{Corollary}

\newcommand{\mP}{\mathbb{P}}
\newcommand{\mE}{\mathbb{E}}

\newcommand{\pr}{\mathbf{P}}
\newcommand{\ev}{\mathbf{E}}

\newcommand{\N}{\mathbb{N}}
\newcommand{\Z}{\mathbb{Z}}
\newcommand{\R}{\mathbb{R}}

\newcommand{\cS}{\mathcal{S}}
\newcommand{\cL}{\mathcal{L}}
\newcommand{\cD}{\mathcal{D}}
\newcommand{\cP}{\mathcal{P}}
\newcommand{\cC}{\mathcal{C}}
\newcommand{\cQ}{\mathcal{Q}}
\newcommand{\cU}{\mathcal{U}}

\newcommand{\di}{\mathrm{d}}
\newcommand{\e}{\mathrm{e}}
\newcommand{\tr}{\operatorname{tr}}

\renewcommand{\mid}{|}

\makeatother

\begin{document}
\begin{frontmatter}

\title{Variational characterization of the critical curve for pinning of random polymers}
\runtitle{The critical curve for pinning of polymers}

\begin{aug}
\author[A]{\fnms{Dimitris} \snm{Cheliotis}\corref{}\ead[label=e1]{dcheliotis@math.uoa.gr}\thanksref{t1}}
\and
\author[B]{\fnms{Frank} \snm{den Hollander}\ead[label=e2]{denholla@math.leidenuniv.nl}}
\thankstext{t1}{Supported in part by the DFG-NWO Bilateral Research
Group ``Mathematical Models from Physics and Biology.''}
\runauthor{D. Cheliotis and F. den Hollander}
\affiliation{University of Athens and Leiden University}
\address[A]{Department of Mathematics\\
University of Athens\\
Panepistimiopolis\\
15784 Athens\\
Greece\\
\printead{e1}} 
\address[B]{Mathematical Institute\\
Leiden University\\
P.O. Box 9512\\
2300 RA Leiden\\
The Netherlands\\
\printead{e2}}
\end{aug}

\received{\smonth{1} \syear{2011}}
\revised{\smonth{7} \syear{2011}}

%
\begin{abstract}
In this paper we look at the pinning of a directed polymer by a one-dimensional
linear interface carrying random charges. There are two phases,
localized and
delocalized, depending on the inverse temperature and on the disorder
bias. Using
quenched and annealed large deviation principles for the empirical
process of
words drawn from a random letter sequence according to a random renewal process
[Birkner, Greven and den Hollander, \textit{Probab. Theory Related Fields}
\textbf{148} (2010) 403--456], we derive variational
formulas for the quenched, respectively, annealed critical curve
separating the
two phases. These variational formulas are used to obtain a necessary
and sufficient
criterion, stated in terms of relative entropies, for the two critical
curves to
be different at a given inverse temperature, a property referred to as relevance
of the disorder. This criterion in turn is used to show that the
regimes of
relevant and irrelevant disorder are separated by a unique inverse critical
temperature. Subsequently, upper and lower bounds are derived for the inverse
critical temperature, from which sufficient conditions under which it
is strictly
positive, respectively, finite are obtained. The former condition is
believed to
be necessary as well, a problem that we will address in a forthcoming paper.

Random pinning has been studied extensively in the literature. The
present paper
opens up a window with a variational view. Our variational formulas for the
quenched and the annealed critical curve are new and provide valuable insight
into the nature of the phase transition. Our results on the inverse critical
temperature drawn from these variational formulas are not new, but they offer
an alternative approach, that is, flexible enough to be extended to
other models
of random polymers with disorder.
\end{abstract}

%
\begin{keyword}[class=AMS]
\kwd[Primary ]{60F10}
\kwd{60K37}
\kwd[; secondary ]{82B27}
\kwd{82B44}.
\end{keyword}
\begin{keyword}
\kwd{Random polymer}
\kwd{random charges}
\kwd{localization vs. delocalization}
\kwd{quenched vs. annealed large deviation principle}
\kwd{quenched vs. annealed critical curve}
\kwd{relevant vs. irrelevant disorder}
\kwd{critical temperature}.
\end{keyword}

\end{frontmatter}
%


\section{Introduction and main results}
\label{S1}


\subsection{Introduction}
\label{S1.1}
\mbox{}

\textit{I. Model}. Let $S=(S_n)_{n\in\N_0}$ be a Markov chain on a
countable state space
$\cS$ in which a given point is marked $0$ ($\N_0=\N\cup\{0\}$). Write
$\pr$ to denote
the law of $S$ given $S_0=0$ and $\ev$ the corresponding expectation.
Let $K$ denote
the distribution of the first return time of $S$ to $0$, that is,
%
\begin{equation}
\label{Kdef}
K(n) := \pr(S_n=0, S_m \neq0\ \forall0<m<n),\qquad n\in\N.
\end{equation}
We will assume that $\sum_{n\in\N} K(n) = 1$ (i.e., $0$ is a
recurrent state)
and
%
\begin{equation}
\label{Kcondgen}
\lim_{n\to\infty} \frac{\log K(n)}{\log n} = -(1+\alpha)
\qquad\mbox{for some } \alpha\in[0,\infty).
\end{equation}

Let $\omega=(\omega_k)_{k\in\N_0}$ be i.i.d. $\R$-valued random
variables with
marginal distribution~$\mu_0$. Write $\mP=\mu_0^{\otimes\N_0}$ to
denote the law
of $\omega$, and $\mE$ to denote the corresponding expectation. We will
assume that
%
\begin{equation}
\label{mgfcond}
M(\lambda) := \mE(\e^{\lambda\omega_0}) < \infty\qquad \forall
\lambda\in\R,
\end{equation}
and that $\mu_0$ has mean $0$ and variance $1$.

Let $\beta\in[0,\infty)$ and $h\in\R$, and for fixed $\omega$ define
the law
$\pr_n^{\beta,h,\omega}$ on $\{0\} \times\cS^n$, the set of $n$-steps
paths in
$\cS$ starting from $0$, by putting
%
\begin{equation}
\label{Pndef}
\frac{\di\pr_n^{\beta,h,\omega}}{\di\pr_n}((S_k)_{k=0}^n)
:= \frac{1}{Z_n^{\beta,h,\omega}} \exp\Biggl[\sum_{k=0}^{n-1}
(\beta\omega_k-h) 1_{\{S_k=0\}}\Biggr] 1_{\{S_n=0\}},
\end{equation}
where $\pr_n$ is the projection of $\pr$ onto $\{0\}\times\cS^n$. Here,
$\beta$ plays
the role of the inverse temperature, $h$ the role of the disorder bias, while
$Z_n^{\beta,h,\omega}$ is the normalizing partition sum. Note that
$k=0$ contributes
to the sum, while $k=n$ does not. Also note that the path is tied to
$0$ at both ends.
This is done for later convenience.

\begin{itremark}
Note that (\ref{Kcondgen}) implies $p:=\mathrm{gcd}[\mathrm
{supp}(K)]=1$. If
$p \geq2$, then the model can be trivially restricted to $p\N$, so
there is no
loss of generality. Moreover, if $\sum_{n\in\N} K(n)<1$, then the model
can be
reduced to the recurrent case by a shift of $h$. Similarly, the
restriction to
$\mu_0$ with mean $0$ and variance $1$ can be removed by a scaling of
$\beta$
and a shift of $h$.
\end{itremark}

\begin{figure}

\includegraphics{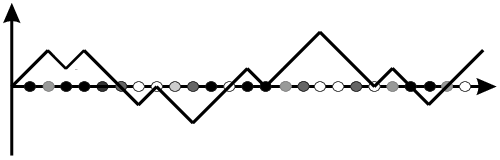}

\caption{A directed polymer sampling random charges at an interface.}
\label{fig-pinning}
\end{figure}

\begin{itremark}
The key example of the above setting is the simple random walk on
$\Z$, for which
$p=2$ and $\alpha=\frac12$ (Spitzer \cite{Sp76}, Section 1). In that
case the
process $(n,S_n)_{n\in\N_0}$ can be thought of as describing a
directed polymer
in $\N_0\times\Z$, that is, pinned to the interface $\N_0 \times\{
0\}$
by random
charges $\beta\omega-h$; see Figure~\ref{fig-pinning}. When the polymer
hits the interface
at time $k$, it picks up a\vadjust{\goodbreak} reward $\exp[\beta\omega_k-h]$, which can be
either $>\!\!1$
or $<\!\!1$, depending on the value of $\omega_k$. For $h \leq0$ the
polymer tends to
intersect the interface with a positive frequency (``localization''),
whereas for
$h>0$ large enough it tends to wander away from the interface
(``delocalization'').
Simple random walk on $\Z^2$ corresponds to $p=2$ and $\alpha=0$,
while simple
random walk on $\Z^d$, $d \geq3$, conditioned on returning to $0$
corresponds to
$p=2$ and $\alpha=\tfrac{d}{2}-1$ (Spitzer \cite{Sp76}, Section 1).
\end{itremark}

\textit{II. Free energy and phase transition}. The \textit{quenched
free energy} is
defined as
%
\begin{equation}
\label{qFreeEnergy}
f^\mathrm{que}(\beta,h) := \lim_{n\to\infty} \frac{1}{n} \log
Z_n^{\beta
,h,\omega}.
\end{equation}
Standard subadditivity arguments show that the limit exists $\omega
$-a.s. and in
$\mP$-mean, and is nonrandom; see, for example, Giacomin~\cite{Gi07},
Chapter 5, and den
Hollander~\cite{dHo09}, Chapter 11. Moreover, $f^\mathrm{que}(\beta,h)
\geq0$ because
$Z_n^{\beta,h,\omega} \geq\e^{\beta\omega_0-h}K(n)$, $n\in\N$, and
$\lim_{n\to\infty}
\frac{1}{n}\log K(n) = 0$ by \eqref{Kcondgen}. The lower bound
$f^\mathrm{que}(\beta,h)
=0$ is attained when $S$ visits the state $0$ only rarely. This
motivates the definition
of two quenched phases,
%
\begin{eqnarray}
\label{LDdefs}
\cL&:=& \{(\beta,h)\dvtx  f^\mathrm{que}(\beta,h)>0\},
\nonumber
\\[-8pt]
\\[-8pt]
\nonumber
\cD&:=& \{(\beta,h)\dvtx  f^\mathrm{que}(\beta,h)=0\},
\end{eqnarray}
referred to as the \textit{localized} phase, respectively, the \textit
{delocalized} phase.

\begin{figure}
\begin{center}
\setlength{\unitlength}{0.35cm}
\begin{picture}(12,12)(0,-1.5)
\put(0,0){\line(12,0){12}}
\put(0,0){\line(0,8){8}}
{\thicklines
\qbezier(0,0)(5,0.5)(9,6.5)
}
\qbezier[60](5,0)(7,3)(11,9)
\put(-.8,-.8){$0$}
\put(12.5,-0.2){$\beta$}
\put(-0.1,8.5){$h$}
\put(0,0){\circle*{.4}}
\put(8,2){$\cL$}
\put(3,3){$\cD$}
\end{picture}
\end{center}
\caption{Qualitative plot of $\beta\mapsto h_c^\mathrm
{que}(\beta)$.
The fine details of this curve are not known.}\label{fig-critcurveranpin}
\end{figure}
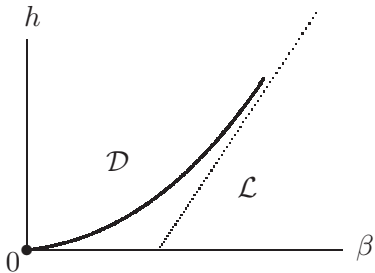

Since $h \mapsto f^\mathrm{que}(\beta,h)$ is nonincreasing for every
$\beta\in
[0,\infty)$, the two phases are separated by a \textit{quenched
critical curve}
%
\begin{equation}
h_c^\mathrm{que}(\beta) := \inf\{h\dvtx  f^\mathrm{que}(\beta
,h)=0\},\qquad
\beta\in[0,\infty),
\end{equation}
with $\cL$ the region below the curve and $\cD$ the region on and
above. Since $(\beta,h)
\mapsto f^\mathrm{que}(\beta,h)$ is convex and $\cD=\{(\beta
,h)\dvtx
f^\mathrm{que}
(\beta,h) \leq0\}$ is a level set of $f^\mathrm{que}$, it follows that
$\cD$ is a
convex set and $h_c^\mathrm{que}$ is a convex function. Since $\beta=0$
corresponds
to a homopolymer, we have $h_c^\mathrm{que}(0)=0$; see \hyperref[A]{Appendix A}. It was
shown in Alexander and Sidoravicius~\cite{AlSi06} that $h_c^\mathrm
{que}(\beta)>0$
for $\beta\in(0,\infty)$. Therefore we have the qualitative picture
drawn in
Figure~\ref{fig-critcurveranpin}. We further remark that $\lim_{\beta
\to
\infty}
h_c^\mathrm{que}(\beta)/\beta$ is finite if and only if $\mathrm
{supp}(\mu_0)$ is
bounded from above.

%
%

The mean value of the disorder is $\mE(\beta\omega_0-h)=-h$. Thus, we
see from
Figure~\ref{fig-critcurveranpin} that for the random pinning model
localization may even
occur for \textit{moderately negative} mean values of the disorder,
contrary to what
happens for the homogeneous pinning model, where localization occurs
only for a strictly
positive parameter; see \hyperref[A]{Appendix~A}. In other words, even a
globally repulsive
random interface can pin the polymer: all that the polymer needs to do
is to hit some
positive values of the disorder and avoid the negative values of the disorder.

The \textit{annealed free energy} is defined by
%
\begin{equation}
f^\mathrm{ann}(\beta,h) := \lim_{n\to\infty} \frac{1}{n}
\log\mE(Z_n^{\beta,h,\omega}).
\end{equation}
Since
%
\begin{equation}
\mE(Z_n^{\beta,h,\omega}) = \ev\Biggl(\exp\Biggl[\sum_{k=0}^{n-1}
[\log M(\beta)-h] 1_{\{S_k=0\}}\Biggr] 1_{\{S_n=0\}}\Biggr),
\end{equation}
we have that $f^\mathrm{ann}(\beta,h)$ is the free energy of the
homopolymer with
parameter $\log M(\beta)-h$. The associated \textit{annealed critical curve}
%
\begin{equation}
h_c^\mathrm{ann}(\beta) := \inf\{h\dvtx  f^\mathrm{ann}(\beta
,h)=0\},\qquad
\beta\in[0,\infty),
\end{equation}
therefore equals
%
\begin{equation}
\label{hcannid}
h_c^\mathrm{ann}(\beta) = \log M(\beta).
\end{equation}

Since $f^\mathrm{que} \leq f^\mathrm{ann}$, we have $h_c^\mathrm{que}
\leq h_c^\mathrm{ann}$.

\begin{itdefinition}
\label{relirreldef}
The disorder is said to be \textit{relevant} for a given choice of
$K$, $\mu_0$
and $\beta$ when $h_c^\mathrm{que}(\beta)<h_c^\mathrm{ann}(\beta)$,
otherwise it is
said to be \textit{irrelevant}.
\end{itdefinition}

\textit{Note}: In the physics literature, the term relevant disorder is
reserved for
the situation where the disorder not only changes the critical value
but also changes
the behavior of the free energy near the critical value. In the present
paper we
adopt the more narrow definition above.\vadjust{\goodbreak}

Our main focus in the present paper will be on deriving \textit{variational
formulas} for $h_c^\mathrm{que}$ and $h_c^\mathrm{ann}$, and on investigating
under what conditions on $K$, $\mu_0$ and $\beta$ the disorder is relevant,
respectively, irrelevant.


\subsection{Main results}
\label{S1.2}

This section contains three theorems and four corollaries, all valid
subject to
(\ref{Kcondgen}) and (\ref{mgfcond}). To state these we need some further notation.

\textit{I. Notation}. Abbreviate
%
\begin{equation}
\label{Edef}
E := \mathrm{supp}[\mu_0] \subset\R.
\end{equation}
Let $\widetilde{E}:= \bigcup_{k\in\N} E^k$ be the set of \textit{finite
words} consisting
of \textit{letters} drawn from~$E$. Let $\cP(\widetilde{E}^\N)$ denote
the set of
probability measures on \textit{infinite sentences},\vspace*{1pt} equipped with the
topology of weak
convergence. Write $\widetilde\theta$ for the left-shift acting on\vspace*{1pt}
$\widetilde{E}^\N$,
and $\cP^{\mathrm{inv}}(\widetilde{E}^\N)$ for the set of probability
measures that are
invariant under $\widetilde{\theta}$.

For $Q\in\cP^\mathrm{inv}(\widetilde{E}^\N)$, let $\pi_{1,1} Q\in
\cP
(E)$ denote the
projection of $Q$ onto the first letter of the first word. Define the set
%
\begin{equation}
\label{Cmu0def}
\cC:= \biggl\{Q\in\cP^\mathrm{inv}(\widetilde E^{\N})\dvtx
\int_E |x| \,\di(\pi_{1, 1} Q)(x) <\infty\biggr\},
\end{equation}
and on this set the function
%
\begin{equation}
\label{PhiDefinition}
\Phi(Q):= \int_E x \,\di(\pi_{1,1} Q)(x),\qquad
Q \in\cC.
\end{equation}
We also need two rate functions on $\cP^\mathrm{inv}(\widetilde
{E}^\N
)$, denoted by
$I^\mathrm{ann}$ and $I^\mathrm{que}$, which will be defined in
Section~\ref{S2}.
These are the rate functions of the annealed and the quenched large
deviation principles
that play a central role in the present paper, and they satisfy
$I^\mathrm{que} \geq
I^\mathrm{ann}$.

\textit{II. Theorems}. With the above ingredients, we obtain the following
\textit{characterization of the critical curves}.

\begin{ittheorem}
\label{qCriticalCurveThm}
Fix $\mu_0$ and $K$. For all $\beta\in[0,\infty)$,
%
\begin{eqnarray}
h_c^\mathrm{que}(\beta)
&=& \sup_{Q \in\cC}[\beta\Phi(Q)-I^\mathrm{que}(Q)],
\label{queCC}\\
h_c^\mathrm{ann}(\beta)
&=& \sup_{Q \in\cC}
[\beta\Phi(Q)-I^\mathrm{ann}(Q)].
\label{annCC}
\end{eqnarray}
\end{ittheorem}

We know that $h_c^\mathrm{ann}(\beta)=\log M(\beta)$. However, the
variational formula for
$h_c^\mathrm{ann}(\beta)$ will be important for the comparison with
$h_c^\mathrm{que}(\beta)$.

Next, for $\beta\in[0,\infty)$ define the probability measures
%
\begin{equation}
\label{mubetadef}
\di\mu_\beta(x) := \frac{1}{M(\beta)} \e^{\beta x} \,\di\mu_0(x),\qquad
x \in E,
\end{equation}
and
%
\begin{eqnarray}
\label{qbetadef}
\di q_\beta(x_1, x_2,\ldots,x_n)
:= K(n)
\,\di\mu_\beta(x_1)\,\di\mu_0(x_2)\times\cdots\times\di\mu_0(x_n),
\nonumber
\\[-8pt]
\\[-8pt]
\eqntext{n\in\N, x_1, x_2,\ldots, x_n\in E.}
\end{eqnarray}
Further, let $Q_\beta:= q_\beta^{\otimes\N} \in\cP^\mathrm
{inv}(\widetilde{E}^\N)$.
Then $Q_0$ is the probability measure under which the words are i.i.d.,
with length
drawn from $K$ and i.i.d. letters drawn from~$\mu_0$, while $Q_\beta$
differs from $Q_0$
in that the first letter of each word is drawn from the tilted
probability distribution
$\mu_\beta$. We will see that $Q_\beta$ is the unique maximizer of the
supremum in
(\ref{annCC}) [note that $Q_\beta\in\cC$ because of (\ref{mgfcond})].
This leads to
the following \textit{necessary and sufficient criterion for disorder
relevance}.

\begin{ittheorem}
\label{CriterionLemma}
Fix $\mu_0$ and $K$. For all $\beta\in[0,\infty)$,
%
\begin{equation}
\label{relirrelcrit}
h_c^\mathrm{que}(\beta) < h_c^\mathrm{ann}(\beta)
\quad\Longleftrightarrow\quad
I^\mathrm{que}(Q_\beta) > I^\mathrm{ann}(Q_\beta).
\end{equation}
\end{ittheorem}

What is appealing about (\ref{relirrelcrit}) is that the gap between
$I^\mathrm{que}$
and $I^\mathrm{ann}$ needs to be established \textit{only} for the
measure $Q_\beta$, which
has a simple and explicit form. We will see that the supremum in (\ref
{queCC}) is attained,
which is to be interpreted as saying that there is a localization
strategy \textit{at}
the quenched critical line.

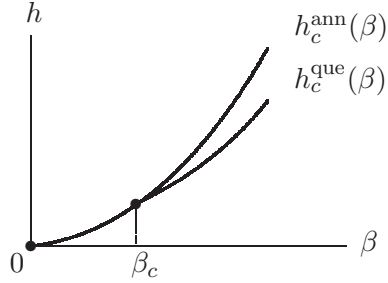
\begin{figure}
\begin{center}
\setlength{\unitlength}{0.35cm}
\begin{picture}(12,12)(0,-1.5)
\put(0,0){\line(12,0){12}}
\put(0,0){\line(0,8){8}}
{\thicklines
\qbezier(0,0)(5,0.5)(9,7.5)
\qbezier(4,1.6)(7,3)(9,5.5)
}
\qbezier[20](4,0)(4,1)(4,1.6)
\put(-.8,-1){$0$}
\put(12.5,-0.2){$\beta$}
\put(-0.2,8.5){$h$}
\put(10,6){$h_c^\mathrm{que}(\beta)$}
\put(10,8){$h_c^\mathrm{ann}(\beta)$}
\put(3.8,-1){$\beta_c$}
\put(0,0){\circle*{.4}}
\put(4,1.6){\circle*{.4}}
\end{picture}
\end{center}
\caption{Uniqueness of the critical inverse temperature $\beta_c$.}
\label{fig-critcurvecomp}
\end{figure}

Disorder relevance is \textit{monotone} in $\beta$; see Figure~\ref
{fig-critcurvecomp}.

\begin{ittheorem}
\label{MonRelevance}
For all $\mu_0$ and $K$ there exists a $\beta_c=\beta_c(\mu_0,K)
\in
[0,\infty]$ such that
%
\begin{equation}
h_c^\mathrm{que}(\beta) \cases{
= h_c^\mathrm{ann}(\beta), & \quad $\mbox{if }
\beta\in[0,\beta_c],$\vspace*{2pt}\cr
< h_c^\mathrm{ann}(\beta), & \quad $\mbox{if } \beta\in(\beta_c,\infty
).$}
\end{equation}
\end{ittheorem}

%


\textit{III. Corollaries}. From Theorems~\ref
{qCriticalCurveThm}--\ref
{MonRelevance}
we draw four corollaries. Abbreviate
%
\begin{equation}\label{chiDefinition}
\chi:= \sum_{n\in\N} [\pr(S_n=0)]^2,\qquad  w := \sup[\mathrm{supp}(\mu_0)].\vadjust{\goodbreak}
\end{equation}

\begin{itcorollary}
\label{disIrrelevancealphazero}
If $\alpha=0$, then $\beta_c=\infty$ for all $\mu_0$.
\end{itcorollary}

\begin{itcorollary}
\label{BdsCritTemp}
If $\alpha\in(0,\infty)$, then the following bounds hold:
\begin{longlist}[(ii)]
\item[(i)] $\beta_c \geq\beta_c^*$ with $\beta_c^* = \beta_c^*(\mu_0,K)
\in
[0,\infty]$
given by
%
\begin{equation}
\label{betac2def}
\beta_c^* := 0\vee\sup\{\beta\dvtx  M(2\beta)/M(\beta)^2 < 1+\chi
^{-1}\}.
\end{equation}
\item[(ii)] $\beta_c \leq\beta_c^{**}$ with $\beta_c^{**}=\beta
_c^{**}(\mu
_0,K) \in
(0,\infty]$ given by
%
\begin{equation} \label{upperbeta}
\beta_c^{**} := \inf\{\beta\dvtx  h(\mu_\beta\mid\mu_0)>h(K)\},
\end{equation}
where $h(\mu_\beta\mid\mu_0)=\int_E \log(\di\mu_\beta/\di\mu
_0) \,\di\mu
_\beta$ is the
relative entropy of $\mu_\beta$ w.r.t. $\mu_0$, and $h(K):=-\sum
_{n\in\N
} K(n)\log K(n)$
is the entropy of $K$.
\end{longlist}
\end{itcorollary}

\begin{itcorollary}
\label{disIrrelevancealphanonzero}
If $\alpha\in(0,\infty)$ and $\chi<\infty$, then $\beta_c>0$ for all
$\mu_0$.
\end{itcorollary}

\begin{itcorollary}
\label{LargeBRelevance}
If $\alpha\in(0,\infty)$, then $\beta_c<\infty$ for all $\mu_0$ with
\mbox{$\mu_0(\{w\})=0$}
(which includes $w=\infty$).
\end{itcorollary}

We close with a conjecture stating that the condition $\chi<\infty$ in
Corollary~\ref{disIrrelevancealphanonzero} is not only sufficient for
$\beta_c>0$ but also necessary. This conjecture will be addressed in
a forthcoming paper.

\begin{itconjecture}
\label{disRelevance}
If $\alpha\in(0,\infty)$ and $\chi=\infty$, then $\beta_c=0$ for all
$\mu_0$.
\end{itconjecture}


\subsection{Discussion}
\label{S1.3}
\mbox{}

\textit{I. What is known from the literature}? Before discussing the
results in
Section~\ref{S1.2}, we give a summary of what is known about the issue
of relevant
vs. irrelevant disorder from the literature. This summary is drawn from
the papers
by Alexander~\cite{Al08}, Toninelli~\cite{To08,To08a},
Giacomin and
Toninelli~\cite{GiTo07}, Derrida, Giacomin, Lacoin and Toninelli~\cite
{DeGiLaTo09},
Alexander and Zygouras~\cite{AlZy09,AlZy10}, Giacomin, Lacoin and
Toninelli~\cite{GiLaTo10,GiLaTo11} and Lacoin~\cite{La}.

\begin{ittheorem}
\label{relirrel}
Suppose that condition (\ref{Kcondgen}) is strengthened to
%
\begin{eqnarray}
\label{Kcond}
K(n) = n^{-(1+\alpha)}L(n)
\nonumber
\\[-8pt]
\\[-8pt]
\eqntext{\mbox{with } \alpha\in[0,\infty)
\mbox{ and } L \mbox{ strictly positive and slowy varying at
infinity}.}
\end{eqnarray}
Then:
\begin{longlist}[(1)]
\item[(1)] $\beta_c=0$ when $\alpha\in(\frac12,\infty)$.

\item[(2)] $\beta_c=0$ when $\alpha=\frac12$ and $\lim_{n\to\infty}
[\log n]^{\delta-1}L^2(n)=0$ for some $\delta>0$.

\item[(3)] $\beta_c>0$ when $\alpha=\frac12$ and
$\sum_{n\in\N} n^{-1} [L(n)]^{-2}<\infty$.

\item[(4)] $\beta_c>0$ when $\alpha\in(0,\frac12)$.

\item[(5)] $\beta_c=\infty$ when $\alpha=0$.\vadjust{\goodbreak}
\end{longlist}
\end{ittheorem}

The results in Theorem~\ref{relirrel} hold irrespective of the choice
of $\mu_0$; see Remark~\ref{disgen} below. Toninelli~\cite{To08a}
proves that if $\log M(\lambda)
\sim C\lambda^\gamma$ as $\lambda\to\infty$ for some $C\in
(0,\infty)$
and $\gamma
\in(1,\infty)$, then $\beta_c<\infty$ irrespective of $\alpha\in
(0,\infty)$ and $L$.
Note that there is a small gap between cases (2) and (3) at the
critical threshold
$\alpha=\frac12$.

For the cases of relevant disorder, bounds on the gap between
$h_c^\mathrm{ann}(\beta)$
and $h_c^\mathrm{que}(\beta)$ have been derived in the above cited
papers subject
to (\ref{Kcond}). As $\beta\downarrow0$, this gap decays like
%
\begin{equation}
\label{gapasymp}
h_c^{\mathrm{ann}}(\beta)-h_c^{\mathrm{que}}(\beta)
\asymp
\cases{
\beta^2, & \quad $\mbox{if } \alpha\in(1,\infty),$\vspace*{2pt}\cr
\beta^2 \psi(1/\beta), & \quad $\mbox{if } \alpha=1,$\vspace*{2pt}\cr
\beta^{2\alpha/(2\alpha-1)}, & \quad $\mbox{if } \alpha\in\bigl(\frac12,1\bigr)$}
\end{equation}
for all choices of $L$, with $\psi$ slowly varying and vanishing at
infinity when
$L(\infty) \in(0,\infty)$.

Partial results are known for $\alpha=\frac12$. For instance, it is
shown in
Giacomin, Lacoin and Toninelli~\cite{GiLaTo11} that, under the
condition in
Theorem~\ref{relirrel}(2), the gap decays faster than any polynomial, namely,
roughly like $\exp[-\beta^{-2/\delta}]$, $\beta\downarrow0$, when $L^2(n)
\asymp[\log n]^{1-\delta}$, $n\to\infty$. This implies that the
disorder can at
most be \textit{marginally relevant}, a situation where standard perturbative
arguments do not work.

\begin{itremark}
\label{disgen}
Some of the above mentioned results are proved for Gaussian
disorder only,
and are claimed to be true for arbitrary disorder subject to (\ref{mgfcond}).
Full proofs for arbitrary disorder are in \cite{DeGiLaTo09,GiLaTo11,La,To08a}.
\end{itremark}

\begin{itremark}
\label{RHarris}
The fact that $\alpha=\frac12$ is critical for relevant vs.
irrelevant disorder
is in accordance with the so-called \textit{Harris criterion} for
disordered systems
(see Harris~\cite{Ha74}): ``Arbitrary weak disorder modifies the nature
of a phase
transition when the order of the phase transition in the nondisordered system
is $<2$.'' The order of the phase transition for the homopolymer, which
is briefly
described in \hyperref[A]{Appendix~A}, is $<2$ precisely when $\alpha\in
(\frac
12,\infty)$
(see Giacomin~\cite{Gi07}, Chapter 2). This link is emphasized in
Toninelli~\cite{To08}.
\end{itremark}

\textit{II. What is new in the present paper}?
The main importance of our results in Section~\ref{S1.2} is that they
\textit{open
up a new window} on the random pinning problem. Whereas the results
cited in
Theorem~\ref{relirrel} are derived with the help of a variety of
\textit{estimation
techniques}, like fractional moment estimates and trial choices of localization
strategies, Theorem~\ref{qCriticalCurveThm} gives a \textit{variational
characterization}
of the critical curves, that is, new. (It is very rare indeed that
critical curves
for disordered systems allow for a direct variational representation.)
Theorem~\ref{CriterionLemma} gives a necessary and sufficient
criterion for
disorder relevance that, although not easy to handle, at least is
explicit and
offers a different handle. Theorem~\ref{MonRelevance} shows that
uniqueness of the
inverse critical\vadjust{\goodbreak} temperature is a direct consequence of this criterion, while
Corollaries~\mbox{\ref{disIrrelevancealphazero}--\ref{LargeBRelevance}} show
that the
criterion can be used to obtain important information on the inverse critical
temperature.

\begin{itremark}
Theorem~\ref{MonRelevance} was proved in Giacomin, Lacoin and
Toni\-nelli~\cite{GiLaTo11} with the help of the FKG-inequality.
\end{itremark}

\begin{itremark}
Corollary~\ref{disIrrelevancealphazero} is the main result in
Alexander and
Zygouras~\cite{AlZy10}.
\end{itremark}

\begin{itremark}
Since (see Section~\ref{S8})
%
\begin{equation}
\label{twolims}
\lim_{\beta\downarrow0} M(2\beta)/M(\beta)^2 = 1,\qquad
\lim_{\beta\to\infty} h(\mu_\beta\mid\mu_0) = \log[1/\mu_0(\{
w\})],
\end{equation}
with the understanding that the second limit is $\infty$ when $\mu
_0(\{
w\})=0$,
Corollary~\ref{BdsCritTemp} implies
Corollaries~\ref{disIrrelevancealphanonzero} and~\ref{LargeBRelevance}.
Corollary~\ref{LargeBRelevance} was noted also in Alexander and
Zygouras~\cite{AlZy10}.
\end{itremark}

\begin{itremark}
\label{I1I2}
Note that $\chi=\ev(|I_1 \cap I_2|)$ with $I_1,I_2$ two
independent copies of
the set of return times of $S$ [recall (\ref{Kdef})]. Thus, according to
Corollary~\ref{disIrrelevancealphanonzero} and Conjecture~\ref{disRelevance},
$\beta_c>0$ is expected to be equivalent to the renewal process of
\textit{joint
return times} to be recurrent. Note that $1/\pr(I_1 \cap I_2 \neq
\varnothing)
= 1+\chi^{-1}$ (see Spitzer~\cite{Sp76}, Section 1), the quantity appearing
in Corollary~\ref{BdsCritTemp}(i).
\end{itremark}

\begin{itremark}
If $\mu_0$ is Bernoulli$(1/2)$ on $\{-1, 1\}$, \eqref{twolims}
gives that
$\lim_{\beta\to\infty} h(\mu_\beta\mid\mu_0)=\log2$.
For any
$\alpha>0$,
we can find a distribution $K$ that satisfies \eqref{Kcondgen} and
$H(K)<\log2$,
and thus \eqref{upperbeta} implies that $\beta_c=\beta_c(\mu_0,
K)<\infty$.
This shows that for $\alpha>0$, the condition $\mu_0(\{w\})=0$ is not
(!) necessary
for $\beta_c<\infty$.
\end{itremark}

\begin{itremark}
As shown in Doney~\cite{Do97}, subject to the condition of regular
variation in
(\ref{Kcond}),
%
\begin{eqnarray}
\pr(S_n=0) \sim\frac{C_\alpha}{n^{1-\alpha}L(n)}
\nonumber
\\[-8pt]
\\[-8pt]
\eqntext{\mbox{as } n\to
\infty
\mbox{ with } C_\alpha=(\alpha/\pi)\sin(\alpha\pi) \mbox{ when
} \alpha
\in(0,1).}
\end{eqnarray}
Hence the condition $\chi<\infty$ in Corollary \ref
{disIrrelevancealphanonzero}
is satisfied exactly for $\alpha\in(0,\frac12)$ and~$L$ arbitrary,
and for $\alpha=\frac12$
and $\sum_{n\in\N} n^{-1}[L(n)]^{-2}<\infty$.\vspace*{2pt} This fits precisely with
cases (3)
and (4) in Theorem~\ref{relirrel}.
\end{itremark}

\begin{itremark}
Corollary~\ref{BdsCritTemp}(ii) is essentially Corollary 3.2 in
Toni\-nelli~\cite{To08a},
where the condition for relevance, $h(\mu_\beta\mid\mu_0)> h(K)$, is
given in an equivalent
form (see equation (3.6) in \cite{To08a}). Note that, by (\ref
{Kcondgen}), $h(K)<\infty$
when $\alpha\in(0,\infty)$.
\end{itremark}


\subsection{Outline}
\label{S1.4}

In Section~\ref{S2} we formulate the annealed and the quenched large
deviation principles
(LDP) that are in Birkner, Greven and den Hollander~\cite{BiGrdHo09},
which are the key
tools in the present paper. In Section~\ref{S3} we use these LDP's to prove
Theorem~\ref{qCriticalCurveThm}. In Section~\ref{S4} we compare the
variational formulas
for the two critical curves and prove the criterion for disorder
relevance stated in
Theorem~\ref{CriterionLemma}. In Section~\ref{S5} we reformulate this
criterion to put
it into a form, that is, more convenient for computations. In
Section~\ref{S6} we use the
latter to prove Theorem~\ref{MonRelevance}. In Sections~\ref
{S7}--\ref
{S8} we prove
Corollaries~\ref{disIrrelevancealphazero}--\ref{LargeBRelevance}.
\hyperref[A]{Appendix~A}
collects a few standard facts about the homopolymer, while
\hyperref[B]{Appendix~B} provides
the details of the proof of a key lemma in Section~\ref{S3} based on an
approximation
argument in \cite{BiGrdHo09}.

\section{Annealed and quenched LDP}
\label{S2}

In this section we recall the main results from Birkner, Greven and den
Hollander~\cite{BiGrdHo09} that are needed in the present paper.
Section~\ref{S2.1}
introduces the relevant notation, while Sections~\ref{S2.2} and \ref
{S2.3} state the
relevant annealed and quenched LDP's.


\begin{figure}

\includegraphics{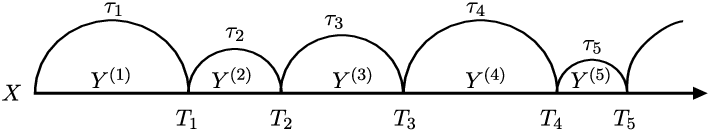}

\caption{Cutting words out from a sequence of letters according to
renewal times.}\label{fig-cutting}
\end{figure}

\subsection{Notation}
\label{S2.1}

Let $E$ be a Polish space, playing the role of an alphabet, that is, a
set of \textit{letters}.
Let $\widetilde{E} := \bigcup_{k\in\N} E^k$ be the set of \textit{finite
words} drawn from $E$,
which can be metrized to become a Polish space.

Fix $\mu_0 \in\cP(E)$, and $K\in\cP(\N)$ satisfying (\ref{Kcondgen}).
Let $X=(X_k)_{k\in\N_0}$
be i.i.d. $E$-valued random variables with marginal law $\mu_0$, and
$\tau=(\tau_i)_{i\in\N}$
i.i.d. $\N$-valued random variables with marginal law $K$. Assume that
$X$ and $\tau$ are
independent, and write $\pr^\ast$ to denote their joint law. Cut words
out of the letter sequence
$X$ according to $\tau$ (see Figure~\ref{fig-cutting}), that is, put
%
\begin{equation}
\label{Tdefs}
T_0:=0\quad \mbox{and}\quad T_i:=T_{i-1}+\tau_i,\qquad i\in\N,
\end{equation}
and let
%
\begin{equation}
\label{eqndefYi}
Y^{(i)} := ( X_{T_{i-1}}, X_{T_{i-1}+1},\ldots, X_{T_{i}-1}),
\qquad i \in\N.
\end{equation}
Under the law $\pr^\ast$, $Y = (Y^{(i)})_{i\in\N}$ is an i.i.d.
sequence of words with marginal
distribution $q_0$ on $\widetilde{E}$ given by
%
\begin{eqnarray}
\label{q0def}
&&\di q_0(x_1,\ldots,x_n)\nonumber\\
&&\qquad := \pr^\ast\bigl(Y^{(1)} \in(\di x_1,\ldots,\di
x_n)\bigr)
\\
&&\qquad\phantom{:}\! = K(n) \,\di\mu_0(x_1) \times\cdots\times\di\mu_0(x_n),\qquad  n\in\N,
x_1,\ldots,x_n\in E.\nonumber
\end{eqnarray}

The reverse operation of \textit{cutting} words out of a sequence of
letters is
\textit{glueing} words together into a sequence of letters.\vspace*{1pt} Formally,
this is done
by defining a \textit{concatenation} map $\kappa$ from $\widetilde
{E}^\N
$ to $E^{\N_0}$.
This map induces in a natural way a map from $\cP(\widetilde{E}^\N)$ to
$\cP(E^{\N_0})$,
the sets of probability measures on $\widetilde{E}^\N$ and $E^{\N_0}$
(endowed with
the topology of weak convergence). The concatenation $q_0^{\otimes\N
}\circ\kappa^{-1}$
of $q_0^{\otimes\N}$ equals $\mu_0^{\N_0}$, as is evident from
(\ref{q0def})


\subsection{Annealed LDP}
\label{S2.2}

Let $\cP^{\mathrm{inv}}(\widetilde{E}^\N)$ be the set of
probability measures
on $\widetilde{E}^\N$ that are invariant under the left-shift
$\widetilde{\theta}$
acting on $\widetilde{E}^\N$. For $N\in\N$, let $(Y^{(1)},\ldots
,Y^{(N)})^\mathrm{per}$
be the periodic extension of the $N$-tuple $(Y^{(1)},\ldots,Y^{(N)})\in
\widetilde{E}^N$
to an element of $\widetilde{E}^\N$, and define
%
\begin{equation}
\label{eqndefRN}
R_N := \frac{1}{N} \sum_{i=0}^{N-1}
\delta_{\widetilde{\theta}^i (Y^{(1)},\ldots,Y^{(N)})^\mathrm{per}}
\in\mathcal{P}^{\mathrm{inv}}(\widetilde{E}^\N).
\end{equation}
This is the \textit{empirical process of $N$-tuples of words}. The
following \textit{annealed LDP}
is standard; see, for example, Dembo and Zeitouni~\cite{DeZe98},
Section 6.5. For $Q\in\cP^\mathrm{inv}
(\widetilde{E}^\N)$, let $H(Q \mid q_0^{\otimes\N})$ be the \textit
{specific relative entropy
of $Q$ w.r.t. $q_0^{\otimes\N}$} defined by
%
\begin{equation}
\label{spentrdef}
H(Q \mid q_0^{\otimes\N}) := \lim_{N\to\infty} \frac{1}{N}
h(\pi_N Q \mid\pi_N q_0^{\otimes\N}),
\end{equation}
where $\pi_N Q \in\cP(\widetilde{E}^N)$ denotes the projection of $Q$
onto the first
$N$ words, $h( \cdot\mid\cdot)$ denotes relative entropy, and the
limit is nondecreasing.

\begin{ittheorem}
\label{aLDP}
The family $\pr^\ast(R_N \in\cdot)$, $N\in\N$, satisfies the
LDP on
$\cP^{\mathrm{inv}}(\widetilde{E}^\N)$ with rate $N$ and with rate function
$I^{\mathrm{ann}}$ given by
%
\begin{equation}
\label{Ianndef}
I^{\mathrm{ann}}(Q):= H(Q \mid q_0^{\otimes\N}),\qquad
Q \in\cP^{\mathrm{inv}}(\widetilde{E}^\N).
\end{equation}
This rate function is lower semi-continuous, has compact level sets,
has a unique
zero at $q_0^{\otimes\N}$, and is affine.
\end{ittheorem}


\subsection{Quenched LDP}
\label{S2.3}

To formulate the quenched analog of Theorem~\ref{aLDP}, we need some
more notation.
Let $\cP^{\mathrm{inv}}(E^{\N_0})$ be the set of probability measures
on $E^{\N_0}$
that are invariant under the left-shift $\theta$ acting on $E^{\N_0}$.
For $Q\in
\cP^{\mathrm{inv}}(\widetilde{E}^\N)$ such that $m_Q := \ev
_{Q}(\tau_1)
< \infty$ (where
$\ev_{Q}$ denotes expectation under the law $Q$ and $\tau_1$ is the
length of the first
word), define
%
\begin{equation}
\label{PsiQdef}
\Psi_Q := \frac{1}{m_Q} \ev_{Q}\Biggl(\sum_{k=0}^{\tau_1-1}
\delta_{\theta^k\kappa(Y)}\Biggr) \in\cP^{\mathrm{inv}}(E^{\N_0}).
\end{equation}
Think of $\Psi_Q$ as the shift-invariant version of $Q\circ\kappa
^{-1}$ obtained
after \textit{randomizing} the location of the origin. This
randomization is necessary
because a shift-invariant~$Q$ in general does not give rise to a shift-invariant
$Q\circ\kappa^{-1}$.

For $\tr\in\N$, let $[\cdot]_{\tr}\dvtx \widetilde{E} \to
[\widetilde
{E}]_{\tr}
= \bigcup_{n=1}^{\tr} E^n$ denote the \textit{truncation map} on words
defined by
%
\begin{equation}
\label{trunword}
\qquad y = (x_1,\ldots,x_n) \mapsto[y]_{\tr}:= (x_1,\ldots,x_{n \wedge\tr}),\qquad
n\in\N, x_1,\ldots,x_n\in E,
\end{equation}
that is, $[y]_{\tr}$ is the word of length $\leq\tr$ obtained from the
word $y$ by dropping
all the letters with label $>\tr$. This map induces in a natural way a
map from
$\widetilde E^\N$ to $[\widetilde{E}]_{\tr}^\N$, and from $\cP
^{\mathrm
{inv}}(\widetilde{E}^\N)$
to $\cP^{\mathrm{inv}}([\widetilde{E}]_{\tr}^\N)$. Note that if
$Q\in\cP
^{\mathrm{inv}}
(\widetilde{E}^\N)$, then $[Q]_{\tr}$ is an element of the set
%
\begin{equation}
\label{Pfin}
\cP^{\mathrm{inv,fin}}(\widetilde{E}^\N) = \{Q\in\cP^{\mathrm
{inv}}(\widetilde{E}^\N)
\dvtx  m_Q<\infty\}.
\end{equation}

\begin{ittheorem}
\label{qLDP}
(Birkner, Greven and den Hollander~\cite{BiGrdHo09})
Assume (\ref{Kcondgen}). Then, for $\mu_0^{\otimes\N_0}$-a.s.
all $X$, the family
of (regular) conditional probability distributions $\mathrm{\pr
}^\ast(R_N \in\cdot\mid X)$,
$N\in\N$, satisfies the LDP on $\cP^{\mathrm{inv}}(\widetilde
{E}^\N)$
with rate $N$ and with
deterministic rate function $I^{\mathrm{que}}$ given by
%
\begin{equation}
\label{eqgndefinitionIalg}
I^\mathrm{que}(Q) := \cases{
I^\mathrm{fin}(Q),
&\quad $\mbox{if } Q
\in\cP^{\mathrm{inv,fin}}(\widetilde{E}^\N),$\vspace*{2pt}\cr
\displaystyle\lim_{\tr\to\infty} I^\mathrm{fin}([Q]_{\tr}),
&\quad $\mbox{otherwise},$}
\end{equation}
where
%
\begin{equation}
\label{eqnratefctexplicitalg}
I^\mathrm{fin}(Q) := H(Q \mid q_0^{\otimes\N})
+ \alpha m_Q H(\Psi_{Q} \mid\mu_0^{\otimes\N_0}).
\end{equation}
This rate function is lower semi-continuous, has compact level sets,
has a unique zero
at $q_0^{\otimes\N}$ and is affine.
\end{ittheorem}

There is no closed form expression for $I^\mathrm{que}(Q)$ when
$m_Q=\infty$. For later
reference we remark that, for all $Q\in\cP^{\mathrm{inv}}(\widetilde
{E}^\N)$,
%
\begin{eqnarray}
\label{truncapproxcont}
I^\mathrm{ann}(Q)
&=& \lim_{\tr\to\infty} I^\mathrm{ann}([Q]_{\tr})
= \sup_{\tr\in\N} I^\mathrm{ann}([Q]_{\tr}),
\nonumber
\\[-8pt]
\\[-8pt]
\nonumber
I^\mathrm{que}(Q)
&=& \lim_{\tr\to\infty} I^\mathrm{que}([Q]_{\tr})
= \sup_{\tr\in\N} I^\mathrm{que}([Q]_{\tr})
\end{eqnarray}
as shown in \cite{BiGrdHo09}, Lemma A.1. A remarkable aspect of (\ref
{eqnratefctexplicitalg}) in
relation to (\ref{Ianndef}) is that it \textit{quantifies} the
difference between $I^{\mathrm{que}}$
and $I^{\mathrm{ann}}$. Note the explicit appearance of the tail
exponent $\alpha$. Also note
that $I^{\mathrm{que}}=I^\mathrm{ann}$ when $\alpha=0$.


\section{\texorpdfstring{Variational formulas: Proof of Theorem~\protect\ref{qCriticalCurveThm}}
{Variational formulas: Proof of Theorem 1.4}}
\label{S3}

In Section~\ref{S3.1} we prove (\ref{annCC}), the variational formula
for the annealed
critical curve. The proof of (\ref{queCC}) in Sections~\ref
{S3.2}--\ref
{S3.4}, the
variational formula for the quenched critical curve, is longer. In
Section~\ref{S3.2}
we first give the proof for $\mu_0$ with finite support. In
Section~\ref
{S3.3} we extend
the proof to $\mu_0$ satisfying (\ref{mgfcond}). In Section~\ref{S3.4}
we prove three
technical lemmas that are needed in Section~\ref{S3.3}.


\subsection{\texorpdfstring{Proof of \protect\eqref{annCC}}{Proof of (1.16)}}
\label{S3.1}
\mbox{}
\begin{pf}
Recall from (\ref{mubetadef}) and (\ref{qbetadef}) that $Q_\beta=q_\beta
^{\otimes\N}$, and
from (\ref{hcannid}) that $h_c^\mathrm{ann}(\beta)=\log M(\beta)$.
Below we show that
for every $Q\in\cP^\mathrm{inv}(\widetilde{E}^\N)$,
%
\begin{equation}
\label{uniqueMinimizer}
\beta\Phi(Q)-I^\mathrm{ann}(Q)=\log M(\beta)-H(Q \mid Q_\beta).
\end{equation}
Taking the supremum over $Q$, we arrive at \eqref{annCC}. Note that the
unique probability
measure that achieves the supremum in (\ref{uniqueMinimizer}) is
$Q_\beta$, which is an
element of the set~$\cC$ defined in (\ref{Cmu0def}) because of (\ref
{mgfcond}).

To get (\ref{uniqueMinimizer}), note that $H(Q \mid Q_\beta)$ is the
limit as $N\to\infty$
of [recall (\ref{mubetadef}) and~(\ref{qbetadef})]
%
\begin{eqnarray}
\label{annentlim}
&&\frac{1}{N} \int_{\widetilde{E}^N} \log\biggl[
\frac{\di(\pi_N Q)}{\di(\pi_N Q_\beta)}
(y_1,\ldots,y_N)\biggr]
\,\di(\pi_N Q) (y_1,\ldots,y_N)\nonumber\\
&&\qquad = \frac{1}{N} \int_{\widetilde{E}^N} \log\biggl[
\frac{\di(\pi_N Q)}{\di(\pi_N Q_0)}
(y_1,\ldots,y_N)\nonumber\\
&&\hspace*{74pt}\qquad{}\times \frac{ M(\beta)^N}{ \e^{\beta[c(y_1)+\cdots+c(y_N)]}}\biggr]
\,\di(\pi_N Q) (y_1,\ldots,y_N)
\\
&&\qquad = \log M(\beta) + \frac{1}{N} h(\pi_NQ \mid\pi_NQ_0)
\nonumber\\
&&\qquad\quad{}- \beta\frac{1}{N} \int_{\widetilde{E}^N}
[c(y_1)+\cdots+c(y_N)] \,\di(\pi_N Q)(y_1,\ldots,y_N),\nonumber
\end{eqnarray}
where, $c(y)$ denotes the first letter of the word $y$. In the last
line of
(\ref{annentlim}), the limit as $N\to\infty$ of the second quantity is
$H(Q \mid Q_0)
=I^\mathrm{ann}(Q)$, while the integral equals $N \Phi(Q)$ by
shift-invariance of $Q$.
Thus, (\ref{uniqueMinimizer}) follows.
\end{pf}


\subsection{\texorpdfstring{Proof of (\protect\ref{queCC}) for $\mu_0$ with finite support}
{Proof of (1.15) for mu0 with finite support}}
\label{S3.2}
\mbox{}
\begin{pf}
The proof comes in three steps.

\textit{Step} 1:
An alternative way to compute the quenched free energy $f^\mathrm
{que}(\beta,h)$ from~(\ref{qFreeEnergy}) is through the radius of convergence $z^\mathrm
{que}(\beta,h)$ of
the power series
%
\begin{equation}
\sum_{n\in\N} z^n Z_n^{\beta,h,\omega},
\end{equation}
because
%
\begin{equation}
\label{zflink}
z^\mathrm{que}(\beta,h) = \e^{-f^\mathrm{que}(\beta,h)}.
\end{equation}
Write
%
\begin{equation}
Z_n^{\beta,h,\omega} = \sum_{N\in\N}
\sum_{0=k_0<k_1<\cdots<k_N=n}
\prod_{i=1}^N K(k_i-k_{i-1}) \e^{\beta\omega_{k_{i-1}}-h},
\end{equation}
so that, for $z\in(0,\infty)$,
%
\begin{equation}
\label{seriesEquality}
\sum_{n\in\N} z^n Z_n^{\beta,h,\omega} = \sum_{N\in\N}
F_N^{\beta
,h,\omega}(z),
\end{equation}
where we abbreviate
%
\begin{equation}
F_N^{\beta,h,\omega}(z) := \sum_{0=k_0<\cdots<k_N<\infty}
\prod_{i=1}^N z ^{k_i-k_{i-1}} K(k_i-k_{i-1}) \e^{\beta\omega_{k_{i-1}}-h}.4
\end{equation}

\textit{Step} 2:
We return to the setting of Section~\ref{S2}. The letter space is $E$,
the word
space is $\widetilde{E} = \bigcup_{k\in\N} E^k$, the sequence of letters
is $\omega
=(\omega_k)_{k\in\N_0}$, while the sequence of renewal times is
$(T_i)_{i\in\N_0}
= (k_i)_{i\in\N_0}$. Each interval $I_i:=[k_{i-1},k_i)$ of integers
cuts out a word
$\omega_{I_i}:=(\omega_{k_{i-1}},\ldots,\omega_{k_i-1})$. Let
%
\begin{equation}
\label{empprocomega}
R_N^\omega= R_N^\omega((k_i)_{i=0}^N)
:= \frac{1}{N} \sum_{i=0}^{N-1}
\delta_{\widetilde{\theta}^i(\omega_{I_1},\ldots,\omega
_{I_N})^{\mathrm{per}}}
\end{equation}
denote the empirical process of $N$-tuples of words in $\omega$ cut out
by the
first $N$ renewals. Then we can rewrite $F_N^{\beta,h,\omega}(z)$ as
%
\begin{eqnarray}
\label{FNid}
F_N^{\beta,h,\omega}(z) &=&
\ev\biggl(\exp\biggl[N\int_{\widetilde{E}}\bigl\{\tau(y)\log z+\bigl(\beta c(y)-h\bigr)\bigr\}
\,\di(\pi_1R_N^\omega)(y)\biggr]\biggr)
\nonumber
\\[-8pt]
\\[-8pt]
\nonumber
&=&\e^{-Nh} \ev\bigl(\exp[Nm_{R_N^\omega}\log z+ N\beta\Phi(R_N^\omega)]\bigr),
\end{eqnarray}
where $\tau(y)$ and $c(y)$ are the length, respectively, the first
letter of the
word $y$, $\pi_1R_N^\omega$ is the projection of $R_N^\omega$ onto the
first word,
while $m_{R_N^\omega}$ and $\Phi(R_N^\omega)$ are the average word
length, respectively,
the average first letter of the first word under~$R_N^\omega$.

To identify the radius of convergence of the series in the left-hand
side of
(\ref{seriesEquality}), we apply the root test for the series in the
right-hand side of
(\ref{seriesEquality}) using the expression in (\ref{FNid}). To that
end, let
%
\begin{equation}
\label{Sdefinition}
S^\mathrm{que}(\beta;z) := \limsup_{N\to\infty} \frac{1}{N} \log
\ev\bigl(\exp[N m_{R_N^\omega} \log z+ N\beta\Phi(R_N^\omega)]\bigr).
\end{equation}
Then
%
\begin{equation}
\label{hSsum}
\limsup_{N\to\infty} \frac{1}{N}\log F_N^{\beta,h,\omega}(z)
= -h + S^\mathrm{que}(\beta;z).
\end{equation}
We know from (\ref{zflink}) and the nonnegativity of $f^\mathrm
{que}(\beta, h)$ that
$z^\mathrm{que}(\beta,h) \leq1$, and we are interested in knowing when
it is $<1$,
respectively, $=1$ [recall (\ref{LDdefs})]. Hence, the sign of the
right-hand side of (\ref{hSsum})
for $z \uparrow1$ will be important as the next lemma shows.

\begin{itlemma}
\label{Scriterion}
For all $\beta\in[0,\infty)$ and $h\in\R$,
%
\begin{eqnarray}
S^\mathrm{que}(\beta; 1-) < h \quad&\Longrightarrow&\quad f(\beta,h)=0,
\nonumber
\\[-8pt]
\\[-8pt]
\nonumber
S^\mathrm{que}(\beta; 1-) > h \quad&\Longrightarrow& \quad f(\beta,h)>0.
\end{eqnarray}
\end{itlemma}

\begin{pf}
The first line holds because, by (\ref{hSsum}), $-h+S^\mathrm
{que}(\beta
;1-)<0$ implies
that the sums in (\ref{seriesEquality}) converge for $|z|<1$, so that
$z^\mathrm{que}
(\beta,h)\geq1$, which gives $f^\mathrm{que}(\beta,h) \leq0$. The
second line holds
because if $-h+S^\mathrm{que}(\beta; 1-)>0$, then there exists a
$z_0<1$ such that
$-h+S^\mathrm{que}(\beta;z_0)>0$, which implies that the sums in
(\ref
{seriesEquality})
diverge for $z=z_0$, so that $z^\mathrm{que}(\beta,h) \leq z_0<1$,
which gives
$f^\mathrm{que}(\beta,h)>0$.
\end{pf}

Lemma~\ref{Scriterion} implies that
%
\begin{equation}
\label{qCritCurve}
h_c^\mathrm{que}(\beta) = S^\mathrm{que}(\beta;1-).
\end{equation}
The rest of the proof is devoted to computing $S^\mathrm{que}(\beta;1-)$.

%
\vspace{1.6cm}
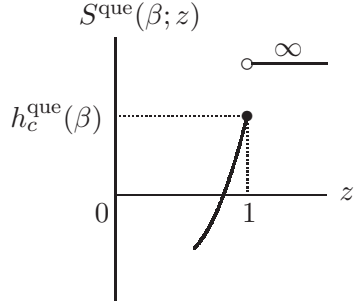
\begin{figure}[t]\vspace*{30pt}
\setlength{\unitlength}{0.35cm}
\begin{picture}(8,6)(0,-1.5)
\put(0,0){\line(8,0){8}}
\put(0,-4){\line(0,10){10}}
{\thicklines
\qbezier(3,-2)(4,-1)(5,3)
\qbezier(5.2,5)(6.5,5)(8,5)
}
\qbezier[30](0,3)(2.5,3)(5,3)
\qbezier[20](5,0)(5,1.5)(5,3)
\put(-.8,-1){$0$}
\put(8.5,-0.2){$z$}
\put(-1.4,6.5){$S^\mathrm{que}(\beta;z)$}
\put(-4,2.6){$h_c^\mathrm{que}(\beta)$}
\put(4.8,-1){$1$}
\put(6,5.3){$\infty$}
\put(5,3){\circle*{.4}}
\put(5,5){\circle{.4}}
\end{picture}
\vspace{0.6cm}
\caption{Qualitative plot of $z \mapsto S^\mathrm{que}(\beta;z)$.}
\label{fig-Scurve}
\end{figure}
\vspace{-1.6cm}

\textit{Step} 3:
Since $\mu_0$ has finite support, $Q\mapsto\Phi(Q)$ is continuous.
Therefore we
can apply Varadhan's lemma to the expression in (\ref{Sdefinition})
for $z=1$
using the LDP of Theorem~\ref{qLDP}. This gives
%
\begin{equation}
\label{Sformulaz=1}
S^\mathrm{que}(\beta;1)
= \sup_{Q\in\cP^\mathrm{inv}(\widetilde{E}^\N)} [\beta\Phi
(Q)-I^\mathrm{que}(Q)].
\end{equation}
We would like to do the same for \eqref{Sdefinition} with $z<1$, and
subsequently
take the limit $z \uparrow1$, to get (see Figure~\ref{fig-Scurve})
%
\begin{equation}
\label{Sformulaz=-1}
S^\mathrm{que}(\beta;1-) = \sup_{Q\in\cP^\mathrm{inv}(\widetilde
{E}^\N)}
[\beta\Phi(Q)-I^\mathrm{que}(Q)].
\end{equation}
However, even though $Q\mapsto\Phi(Q)$ is continuous (because $\mu_0$
has finite
support), $Q \mapsto m_Q$ is only lower semicontinuous. Therefore we
proceed by
first showing that the term $Nm_{R_N^\omega}\log z$ in (\ref
{Sdefinition}) is
harmless in the limit as $z \uparrow1$.

\begin{itlemma}
\label{lsc_lemma}
$S^\mathrm{que}(\beta;1-)=S^\mathrm{que}(\beta;1)$ for all $\beta
\in
[0,\infty)$.
\end{itlemma}

\begin{pf}
Since $S^\mathrm{que}(\beta;1-) \leq S^\mathrm{que}(\beta;1)$, we need
only prove the
reverse inequality. The idea is to show that, for any $Q\in\cP
^\mathrm
{inv}(\widetilde{E}^\N)$
and in the limit as $N\to\infty$, $R_N^\omega$ can be arbitrarily close
to $Q$ with
probability $\approx\exp[-NI^\mathrm{que}(Q)]$ while $m_{R_N^\omega}$
remains bounded
by a large constant. Therefore, letting $N\to\infty$ followed by
$z\uparrow1$, we can
remove the term $Nm_{R_N^\omega}\log z$ in (\ref{Sdefinition}). The
details are given in
\hyperref[B]{Appendix~B}.
\end{pf}

Combining Lemma~\ref{lsc_lemma} with \eqref{qCritCurve} and \eqref
{Sformulaz=1}, we
obtain \eqref{queCC}.
\end{pf}


\subsection{\texorpdfstring{Proof of (\protect\ref{queCC}) for $\mu_0$ satisfying (\protect\ref{mgfcond})}
{Proof of (1.15) for mu0 satisfying (1.3)}}
\label{S3.3}

The proof stays the same up to~(\ref{qCritCurve}). Henceforth write
$\cC
=\cC(\mu_0)$
to exhibit the fact that the set $\cC$ in (\ref{Cmu0def}) depends on
$\mu_0$ via
its support $E$ in (\ref{Edef}), and define
%
\begin{equation}
\label{Adef}
A(\beta) := \sup_{Q\in\cC(\mu_0)} [\beta\Phi(Q)-I^\mathrm{que}(Q)],
\end{equation}
which replaces the right-hand side of (\ref{Sformulaz=-1}). We will
show the following.

\begin{itlemma}
\label{leftLimit}
$S^\mathrm{que}(\beta;1-)=A(\beta)$ for all $\beta\in(0,\infty)$.
\end{itlemma}

\begin{pf}
The proof of the lemma is accomplished in four steps. Along the way we
use three technical
lemmas, the proof of which is deferred to Section~\ref{S3.4}. Our
starting point is the
validity of the claim for $\mu_0$ with finite support obtained in
Lemma~\ref{lsc_lemma}.
(Note that $|E|<\infty$ implies $\cC=\cC(\mu_0)=\cP^\mathrm
{inv}(\widetilde E^{\N})$.)

\textit{Step} 1: $S^\mathrm{que}(\beta;1-) \leq A(\beta)$ for all
$\beta\in(0,\infty)$
when $\mu_0$ satisfies (\ref{mgfcond}).

\begin{pf}
We have $S^\mathrm{que}(\beta;1-) \leq S^\mathrm{que}(\beta;1)$. We
will show that
$S^\mathrm{que}(\beta;1) \leq A(p\beta)/p$ for all $p>1$. Taking $p
\downarrow1$ and
using the continuity of $A$, proven in Lemma~\ref{ACont} below, we get
the claim.

For $M>0$, let
%
\begin{equation}
\Phi^M(Q) := \int_E (x\wedge M) \,\di(\pi_{1, 1} Q)(x).
\end{equation}
Then, for any $p,q>1$ such that $p^{-1}+q^{-1}=1$, we have
%
\begin{eqnarray}
\ev\bigl(\e^{N \beta\Phi(R_N^\omega)}\bigr)
&=& \ev\bigl(\e^{\beta\sum_{i=1}^N c(y_i) 1_{\{c(y_i)\leq M\}}}
\e^{\beta\sum_{i=1}^N c(y_i) 1_{\{c(y_i)>M\}}}\bigr)\nonumber\\
&\leq&\bigl[\ev\bigl(\e^{p\beta\sum_{i=1}^N c(y_i) 1_{\{c(y_i)\leq M\}}}\bigr)\bigr]^{1/p}
\bigl[\ev\bigl(\e^{q\beta\sum_{i=1}^N c(y_i) 1_{\{c(y_i)> M\}}}\bigr)\bigr]^{1/q}\\
&\leq&\bigl[\ev\bigl(\e^{Np\beta\Phi^M(R_N^\omega)}\bigr)\bigr]^{1/p}
\bigl[\ev\bigl(\e^{q\beta\sum_{i=1}^N c(y_i) 1_{\{c(y_i)>
M\}}}\bigr)\bigr]^{1/q},\nonumber
\end{eqnarray}
where $y_1,\ldots,y_N$ are the $N$ words determining $R_N^\omega$ and
$c(y_i)$ is the first
letter of the $i$th word. Hence
%
\begin{eqnarray}
\label{truncationUp}
\frac{1}{N} \log\ev\bigl(\e^{N\beta\Phi(R_N^\omega)}\bigr)
&\leq&\frac{1}{p} \frac{1}{N} \log\ev\bigl(\e^{Np\beta\Phi
^M(R_N^\omega)}\bigr)
\nonumber
\\[-8pt]
\\[-8pt]
\nonumber
&&{}+ \frac{1}{q} \frac{1}{N} \log\ev\bigl(\e^{q\beta\sum_{i=1}^N c(y_i)
1_{\{c(y_i)>M\}}}\bigr).
\end{eqnarray}
Since $Q \mapsto\Phi^M(Q)$ is upper semicontinuous, Varadhan's lemma gives
%
\begin{equation}
\label{truncationUpext}
\quad\limsup_{N\to\infty} \frac{1}{N} \log\ev\bigl(\e^{Np\beta\Phi
^M(R_N^\omega)}\bigr)
\leq\sup_{Q\in\cP^\mathrm{inv}(\widetilde E^{\N})} [p \beta\Phi
^M(Q)-I^\mathrm{que}(Q)].
\end{equation}
Clearly, $Q$'s with $\int_E (x\wedge0) \,\di(\pi_{1, 1} Q)(x) =
-\infty$
do not contribute
to the supremum. Also, $Q$'s with $\int_E (x \vee0) \,\di(\pi_{1,1}
Q)(x) = \infty$ do not
contribute, because for such $Q$ we have $I^\mathrm{que}(Q)=\infty$, by
Lemma~\ref{FiniteEntropy}
below, and $\Phi^M(Q)<\infty$. Since $\Phi^M \leq\Phi$, we
therefore have
%
\begin{eqnarray}
\sup_{Q\in\cP^\mathrm{inv}(\widetilde E^{\N})} [p \beta\Phi
^M(Q)-I^\mathrm{que}(Q)]
&\leq&\sup_{Q\in\cC(\mu_0)} [p \beta\Phi(Q)-I^\mathrm{que}(Q)]
\nonumber
\\[-8pt]
\\[-8pt]
\nonumber
& =&
A(p\beta).
\end{eqnarray}

Next, we use the following observation. For any sequence $\Theta
=(\Theta
_N)_{N\in\N}$
of positive random variables on a space with probability measure $\mP$,
we have
%
\begin{equation}
\label{annealedBound}
\limsup_{N\to\infty} \frac{1}{N} \log\Theta_N \leq
\limsup_{N\to\infty} \frac{1}{N} \log\mE(\Theta_N)\qquad \mP\mbox{-a.s.},
\end{equation}
by the first Borel--Cantelli lemma. Applying this to
%
\begin{eqnarray}
\Theta_N := \ev\bigl(\e^{q\beta\sum_{i=1}^N c(y_i)
1_{\{c(y_i)>M\}}}\bigr)\\
\eqntext{\mbox{with }
\mE(\Theta_N) = \biggl(\displaystyle\int_E
\e^{q\beta x 1_{\{x>M\}}} \,\di\mu_0(x)\biggr)^N =:(c_M)^N,}
\end{eqnarray}
we get, after letting $N\to\infty$ in \eqref{truncationUp},
%
\begin{equation}
\label{resann}
S^\mathrm{que}(\beta;1) \leq\frac{1}{p} A(p\beta) + \frac{1}{q}
\log c_M.
\end{equation}
By (\ref{mgfcond}), we have $c_M<\infty$ for all $M>0$ and $\lim
_{M\to
\infty}c_M=1$.
Hence\break $S^\mathrm{que}(\beta;1) \leq A(p\beta)/p$.
\end{pf}

\textit{Step} 2: $S^\mathrm{que}(\beta;1-) \geq A(\beta)$ for all
$\beta\in(0,\infty)$
when $\mu_0$ has bounded support.

\begin{pf}
In the estimates below, we abbreviate
%
\begin{equation}
L^\omega_N := N m_{R_N^\omega},
\end{equation}
the sum of the lengths of the first $N$ words. The proof is based on a
discretization
argument similar to the one used in \cite{BiGrdHo09}, Section 8. For
$\delta>0$ and
$x\in E$, let $\langle x\rangle_\delta:=\sup\{k\delta\dvtx  k\in\Z,k
\delta\leq x\}$.
The operation $\langle\cdot\rangle$ extends to measures on $E$,
$\widetilde E$ and
$\widetilde E^\N$ in the obvious way. Now, $\langle R_N^\omega\rangle
_\delta$ satisfies
the quenched LDP with rate function $I^\mathrm{que}_\delta$, the
quenched rate function
corresponding to the measure $\langle\mu_0\rangle_\delta$. Clearly,
%
\begin{equation}
\ev\bigl(\e^{L^\omega_N\log z+N\beta\Phi(R_N^\omega)}\bigr)
\geq\ev\bigl(\e^{L^\omega_N\log z+N\beta\Phi(\langle R_N^\omega
\rangle
_\delta)}\bigr),
\end{equation}
and so, by the results in Section~\ref{S3.2}, we have
%
\begin{equation}
\label{Slbinfsup}
S^\mathrm{que}(\beta; 1-) \geq\sup_{Q \in\cC(\langle\mu
_0\rangle
_\delta)}
[\beta\Phi(Q)-I^\mathrm{que}_\delta(Q)].
\end{equation}
For every $Q\in\cC(\mu_0)$, we have
%
\begin{equation}
\Phi(Q) = \lim_{\delta\downarrow0} \Phi(\langle Q\rangle_\delta),
\qquad
I^\mathrm{que}(Q) = \lim_{n\to\infty}
I^\mathrm{que}_{\delta_n}(\langle Q\rangle_{\delta_n}),
\end{equation}
where $\delta_n=2^{-n}$. The first relation holds because $\Phi
(\langle
Q\rangle_\delta) \leq
\Phi(Q) \leq\Phi(\langle Q\rangle_\delta)+\delta$, the second
relation uses
Lemma~\ref{rateFunctionApprox}(i) below. Hence the claim follows by
picking $\delta=\delta_n$
in (\ref{Slbinfsup}) and letting $n\to\infty$.
\end{pf}

\textit{Step} 3: $S^\mathrm{que}(\beta;1-) \geq A(\beta)$ for all
$\beta\in(0,\infty)$
when $\mu_0$ satisfies (\ref{mgfcond}) with support bounded from below.

\begin{pf}
For $M>0$ and $x\in E$, let $x^M = x \wedge M$. This truncation
operation acts on $\mu_0$
by moving the mass in $(M,\infty)$ to $M$, resulting in a measure $\mu
_0^M$ with bounded
support and with associated quenched rate function $I^{\mathrm
{que},M}$. Let $R_N^{\omega,M}$
be the empirical process of $N$-tuples of words obtained from
$R_N^\omega$ defined in
(\ref{eqndefRN}) after replacing each letter $x\in E$ by $x^M$. We have
%
\begin{equation}
\ev\bigl(e^{L^\omega_N\log z+N\beta\Phi(R_N^\omega)}\bigr)
\geq
\ev\bigl(\e^{L^\omega_N\log z+N\beta\Phi(R_N^{\omega, M})}\bigr).
\end{equation}
Combined with the result in Step 2, this bound implies that
%
\begin{equation}
\label{Slbinfsupp}
S(\beta; 1-) \geq\sup_{Q'\in\cC(\mu_0^M)} [\beta\Phi
(Q')-I^{\mathrm
{que},M}(Q')].
\end{equation}
For every $Q\in\cC(\mu_0)$, we have
%
\begin{eqnarray}
\label{PhiIlims}
\Phi(Q) &=& \lim_{M\to\infty} \Phi(Q^M)
= \lim_{M\to\infty} \int_E (x\wedge M) \,\di(\pi_{1,1} Q)(x),
\nonumber
\\[-8pt]
\\[-8pt]
\nonumber
I^\mathrm{que}(Q) &=& \lim_{M\to\infty} I^{\mathrm{que},M}(Q^M).
\end{eqnarray}
The first relation holds by dominated convergence, and the second
relation uses
Lemma~\ref{rateFunctionApprox}(ii) below. It follows from (\ref{PhiIlims})
that
%
\begin{eqnarray}
\limsup_{M\to\infty} \sup_{Q'\in\cC(\mu_0^M)} [\beta\Phi
(Q')-I^{\mathrm{que},M}(Q')]
\geq\beta\Phi(Q)-I^{\mathrm{que}}(Q)
\nonumber
\\[-8pt]
\\[-8pt]
\eqntext{\forall Q\in\cC(\mu_0),}
\end{eqnarray}
which combined with (\ref{Slbinfsupp}) yields
%
\begin{equation}
S(\beta; 1-) \geq\beta\Phi(Q)-I^{\mathrm{que}}(Q) \qquad\forall Q\in\cC
(\mu_0).
\end{equation}
Take the supremum over $Q\in\cC(\mu_0)$ to get the claim.
\end{pf}

\textit{Step} 4: $S^\mathrm{que}(\beta;1-) \geq A(\beta)$ for all
$\beta\in(0,\infty)$
when $\mu_0$ satisfies (\ref{mgfcond}).

\begin{pf}
For $M>0$ and $x\in E$, let $x^{-M}=x\vee(-M)$. This truncation
operation acts on $\mu_0$
by moving the mass in $(-\infty,-M)$ to $-M$, resulting in a measure
$\mu_0^{-M}$ with
support bounded from below and with associated quenched rate function
$I^{\mathrm{que},-M}$.
Let $R_N^{\omega,-M}$ be the empirical process of $N$-tuples of words
obtained from
$R_N^\omega$ defined in (\ref{eqndefRN}) after replacing each letter
$x\in E$ by
$x^{-M}$.

As in Step 1, for any $p,q>1$ such that $p^{-1}+q^{-1}=1$, we have
%
\begin{eqnarray}
&&\ev\bigl(\e^{L^\omega_N\log z+N\beta\Phi(R_N^{\omega, -M})}\bigr)\nonumber\\
&&\qquad\leq\ev\bigl(\e^{L^\omega_N\log z+N\beta\Phi(R_N^\omega)} \e
^{-\beta\sum_{i=1}^N
c(y_i) 1_{\{c(y_i)<-M\}}}\bigr)\\
&&\qquad\leq\bigl[\ev\bigl(\e^{pL^\omega_N\log z + Np\beta\Phi(R_N^\omega)}\bigr)\bigr]^{1/p}
\bigl[\ev\bigl(\e^{-q\beta\sum_{i=1}^N c(y_i)
1_{\{c(y_i)<-M\}}}\bigr)\bigr]^{1/q},\nonumber
\end{eqnarray}
and hence
%
\begin{eqnarray}
&&\frac{1}{N} \log\ev\bigl(\e^{L^\omega_N\log z+N\beta\Phi(R_N^{\omega,
-M})}\bigr)\nonumber\\
&&\qquad \leq\frac{1}{p} \frac{1}{N}
\log\ev\bigl(\e^{p L^\omega_N\log z+Np\beta\Phi(R_N^\omega)}\bigr)
\\
&&\qquad\quad{}+ \frac{1}{q} \frac{1}{N}
\log\ev\bigl(\e^{-q\beta\sum_{i=1}^N c(y_i)
1_{\{c(y_i)<-M\}}}\bigr).\nonumber
\end{eqnarray}
Let $N\to\infty$ followed by $z \uparrow1$. For the left-hand side, we
have the lower
bound in Step 3, while the second term in the right-hand side can be
handled as in
(\ref{annealedBound}--\ref{resann}). Therefore, recalling (\ref
{Sdefinition}) and
writing $p\log z=\log z^p$, we get
%
\begin{eqnarray}
\sup_{Q\in\cC(\mu_0^{-M})} [\beta\Phi(Q)-I^{\mathrm{que}, -M}(Q)]
\leq\frac{1}{p} S^\mathrm{que}(p \beta; 1-) + \frac{1}{q} \log
C_{-M}
\nonumber
\\[-8pt]
\\[-8pt]
\eqntext{\mbox{ with } C_{-M} := \displaystyle\int_E \mathrm{e}^{-q\beta x 1_{\{x<-M\}}}
\,\di\mu_0(x).}
\end{eqnarray}
Letting $M\to\infty$ and using that $\lim_{M\to\infty} C_{-M}=1$ by
(\ref{mgfcond}),
we arrive at
%
\begin{equation}
\qquad\frac{1}{p} S^\mathrm{que}(p\beta, 1-) \geq\limsup_{M\to\infty}
\sup_{Q\in\cC(\mu_0^{-M})} [\beta\Phi(Q)-I^{\mathrm{que}, -M}(Q)]
\geq A(\beta),
\end{equation}
where the last inequality is obtained via arguments similar to those following~(\ref{Slbinfsupp}), which require the use of Lemma~\ref
{rateFunctionApprox}(iii)
below. Finally, let $p \downarrow1$, and use the continuity of $\beta
\mapsto S(\beta; 1-)$,
proven in Lemma~\ref{ACont} below.
\end{pf}

This completes the proof of Lemma~\ref{leftLimit} and hence of
Theorem~\ref{qCriticalCurveThm}.
\end{pf}


\subsection{Technical lemmas}
\label{S3.4}

In the proof of Lemma~\ref{leftLimit} we used three technical lemmas,
which we prove
in this section.

\begin{itlemma}
\label{ACont}
$\beta\mapsto A(\beta)$ and $\beta\mapsto S^\mathrm{que}(\beta
;1-)$ are finite
and convex on $[0,\infty)$ and, consequently, are continuous on
$(0,\infty)$.
\end{itlemma}

\begin{pf}
For the first function, note that $A(\beta) \leq\sup_{Q\in\cC(\mu_0)}
[\beta\Phi(Q)-
I^\mathrm{ann}(Q)] \leq\log M(\beta)<\infty$ by (\ref{mgfcond}) and
(\ref{uniqueMinimizer}),
and convexity follows from the fact that $A$ is a supremum of linear
functions. For the
second function, note that $S^\mathrm{que}(\beta;1-) \leq S^\mathrm
{que}(\beta;1) = A(\beta)$,
and convexity follows from H\"older's inequality.
\end{pf}

\begin{itlemma}
\label{FiniteEntropy}
If $\mu,\nu\in P(\R)$ satisfy $h(\mu\mid\nu)<\infty$ and $\int
_E \e
^{\lambda x}
\,\mathrm{d}\nu(x)<\infty$ for some $\lambda>0$, then $\int_E (x
\vee0)
\,\mathrm{d}
\mu(x)<\infty$.
\end{itlemma}

\begin{pf}
The claim follows from the inequality
%
\begin{equation}
\label{LegendreTransIneq}
\int_E f \,\di\mu\leq h(\mu\mid\nu) + \log\int_E \e^f \,\di\nu,
\end{equation}
which is valid for all bounded and measurable $f$ (see Dembo and
Zeitouni~\cite{DeZe98},
Lemma 6.2.13) and, by monotone convergence, extends to measurable $f
\geq0$. Pick
$f(x) = \lambda(x \vee0)$, $x\in E$.
\end{pf}

\begin{itlemma}
\label{rateFunctionApprox}
For every $Q\in\cP^\mathrm{inv}(\widetilde{E}^\N)$:
\begin{longlist}[(iii)]
\item[(i)] $\lim_{n\to\infty} I^\mathrm{que}_{\delta_n}(\langle Q\rangle
_{\delta_n})= I^\mathrm{que}(Q)$ with $\delta_n:=2^{-n}$;

\item[(ii)] $\lim_{M\to\infty} I^{\mathrm{que},M}(Q^M)=I^\mathrm
{que}(Q)$;

\item[(iii)] $\lim_{M\to\infty} I^{\mathrm{que}, -M}(Q^{-M})=I^\mathrm{que}(Q)$.
\end{longlist}
\end{itlemma}

\begin{pf}
(i) The proof proceeds by choosing an appropriate function $I \dvtx
[0,1]\to\R$
and proving that:
%
\begin{eqnarray}
\label{alpbet}
\mathrm{(a)}&& \quad I(0)=\lim_{\delta\downarrow0} I(\delta);
\nonumber\\
\mathrm{(b)} &&\quad I(0)\geq I(\delta_1) \geq I(\delta_2)\hspace*{140pt}\\
\eqntext{\mbox{whenever }
\delta_2=k \delta_1 \in(0,1)
\mbox{ for some } k\in\N.}
\end{eqnarray}
Recalling (\ref{eqgndefinitionIalg}) and (\ref{eqnratefctexplicitalg}), we
see that we need
the following choices for~$I$:
%
\begin{eqnarray}
\label{alpbetext}
\mathrm{(1)}\quad I(\delta) &=& \cases{
N^{-1}h(\langle\pi_N Q\rangle_\delta\mid\langle\pi_N
q_0^{\otimes\N
}\rangle_\delta),
&\quad $\delta>0,$\vspace*{2pt}\cr
N^{-1}h(\pi_N Q \mid\pi_N q_0^{\otimes\N}),
&\quad $\delta=0,$}
\nonumber\\
\mathrm{(2)}\quad I(\delta) &=& \cases{
H(\langle Q \rangle_\delta\mid\langle q_0^{\otimes\N}\rangle
_\delta),
&\quad $\delta>0,$\vspace*{2pt}\cr
H(Q \mid q_0^{\otimes\N}),
&\quad $\delta=0,$}
\nonumber
\\[-8pt]
\\[-8pt]
\nonumber
\mathrm{(3)} \quad I(\delta) &=& \cases{
N^{-1}h(\langle\pi_N\Psi_Q\rangle_\delta\mid
\langle\pi_N \mu_0^{\otimes\N_0} \rangle_\delta),
&\quad $\delta>0,$\vspace*{2pt}\cr
N^{-1}h(\pi_N\Psi_Q \mid\pi_N\mu_0^{\otimes\N_0}),
&\quad $\delta=0,$}\\
\mathrm{(4)}\quad I(\delta) &=& \cases{
H(\langle\Psi_Q \rangle_\delta\mid\langle\mu_0^{\otimes\N
_0}\rangle
_\delta),
&\quad $\delta>0,$\vspace*{2pt}\cr
H(\Psi_Q \mid\mu_0^{\otimes\N_0} ),
&\quad $\delta=0,$}\nonumber
\end{eqnarray}
with $N\in\mathbb{N}$.
It is clear from the definition of specific relative entropy [recall~\ref{spentrdef})]
that if (a) and (b) hold for the choices (1) and (3), then they also
hold for the
choices~(2) and (4), respectively. We will not actually prove (a) and
(b) for the
choices (1) and (3), but for the simpler choice
%
\begin{equation}
\label{simch}
I(\delta) = \cases{
h(\langle\mu\rangle_\delta\mid\langle\mu_0 \rangle_\delta),
&\quad $\delta>0,$\vspace*{2pt}\cr
h(\mu\mid\mu_0),
&\quad $\delta=0.$}
\end{equation}
The proof will make it evident how to properly deal with (1) and (3).

Let $B(\R)$ be the set of real-valued, bounded and Borel measurable
functions on~$\R$ and,
for $\phi\in B(\R)$ and $\delta>0$, let $\phi_\delta$ be the function
defined by
$\phi_\delta(x):= \phi(\langle x \rangle_\delta)$. As shown in
Dembo and
Zeitouni~\cite{DeZe98}, Lemma 6.2.13, we have
%
\begin{eqnarray}
\label{hsuprep}
h(\langle\mu\rangle_\delta\mid\langle\mu_0 \rangle_\delta)
&= &\sup_{\phi\in B(\R)} \biggl\{\int_\R\phi\,\di\langle\mu\rangle
_\delta
-\log\int_\R\e^\phi\,\di\langle\mu_0 \rangle_\delta\biggr\}
\nonumber
\\[-8pt]
\\[-8pt]
\nonumber
&=& \sup_{\phi\in B(\R)} \biggl\{\int_\R\phi_\delta\,\di\mu
-\log\int_\R\e^{\phi_\delta} \,\di\mu_0\biggr\}.
\end{eqnarray}
From this representation, property (b) follows for the choice in (\ref
{simch}). Next,
fix any $\varepsilon>0$ and take a $\phi$ such that $\int_\R\phi
\,\di\mu
-\log\int_\R
\e^\phi\,\di\mu_0 \geq h(\mu\mid\mu_0)-\varepsilon$. Then, since~$\phi
_\delta$ converges
pointwise to $\phi$ as $\delta\downarrow0$, the bounded convergence
theorem together
with (\ref{hsuprep}) give
%
\begin{equation}
\liminf_{\delta\downarrow0} h(\langle\mu\rangle_\delta\mid
\langle\mu_0 \rangle_\delta) \geq h(\mu\mid\mu_0)-\varepsilon.
\end{equation}
Hence $\liminf_{\delta\downarrow0} I(\delta) \geq I(0)-\varepsilon$.
Since $I(0) \geq
I(\delta)$, property (a) follows after letting $\varepsilon\downarrow0$.

Having thus convinced ourselves that (\ref{alpbet}) and (\ref{alpbetext})
are true, we now
know that for any $Q\in\cP^\mathrm{inv}(\widetilde{E}^\N)$ the sequences
%
\begin{equation}
H(\langle Q \rangle_{\delta_n} \mid\langle q_0^{\otimes\N}\rangle
_{\delta_n}),
\qquad
H(\langle\Psi_Q \rangle_{\delta_n} \mid\langle\mu_0^{\otimes\N
_0}\rangle_{\delta_n}),\qquad
n\in\N,
\end{equation}
are increasing and converge to $H(Q \mid q_0^{\otimes\N})$, respectively,
$H(\Psi_Q \mid\mu_0^{\otimes\N_0})$. This implies the claim for $Q$
with $m_Q<\infty$ [recall
(\ref{eqnratefctexplicitalg})].\vspace*{1pt} For $Q$ with $m_Q=\infty$ we use that
$I^\mathrm{que}(Q)
= \sup_{\tr\in\N} I([Q]_{\tr})$ [recall (\ref{truncapproxcont})], to
conclude that
$I^\mathrm{que}_{\delta_n}(\langle Q\rangle_{\delta_n})$ is increasing
and converges
to $I^\mathrm{que}(Q)$.

(ii)--(iii) The proof is similar as for (i).
\end{pf}


\section{\texorpdfstring{Characterization of disorder relevance: Proof of Theorem~\protect\ref{CriterionLemma}}
{Characterization of disorder relevance: Proof of Theorem 1.5}}\label{S4}
\mbox{}
\begin{pf}
We will need the following lemma, the proof of which is postponed.

\begin{itlemma}
\label{supAttained}
The supremum $\sup_{Q\in\cC}[\beta\Phi(Q)-I^\mathrm{que}(Q)]$ is
attained for all
$\beta\in(0,\infty)$.
\end{itlemma}

Let $Q^*$ be a measure achieving the supremum in Lemma \ref
{supAttained}. Suppose that
$h_c^\mathrm{que}(\beta) = h_c^\mathrm{ann}(\beta)$. Then
%
\begin{eqnarray}
\label{Icomps}
h_c^\mathrm{que}(\beta) &= &\beta\Phi(Q^*)-I^\mathrm{que}(Q^*)
\leq\beta\Phi(Q^*)-I^\mathrm{ann}(Q^*)
\nonumber
\\[-8pt]
\\[-8pt]
\nonumber
& \leq&\beta\Phi(Q_\beta)-I^\mathrm{ann}(Q_\beta)
= h_c^\mathrm{ann}(\beta) = h_c^\mathrm{que}(\beta),
\end{eqnarray}
where the second equality uses that $Q_\beta$ achieves the supremum in
\eqref{annCC}
[with $I^\mathrm{ann}(Q_\beta)<\infty$], as shown by \eqref
{uniqueMinimizer}. It follows
that both inequalities in (\ref{Icomps}) are equalities. However, since
$Q_\beta$ uniquely
achieves the supremum in \eqref{annCC}, we must have $Q^*=Q_\beta$
and therefore
$I^\mathrm{que}(Q_\beta) = I^\mathrm{ann} (Q_\beta)$.

Conversely, suppose that $I^\mathrm{que}(Q_\beta)=I^\mathrm
{ann}(Q_\beta
)$. Then
%
\begin{equation}
\qquad h_c^\mathrm{que}(\beta) \geq[\beta\Phi(Q_\beta) - I^\mathrm
{que}(Q_\beta)]
= [\beta\Phi(Q_\beta)- I^\mathrm{ann}(Q_\beta)] = h_c^\mathrm
{ann}(\beta).
\end{equation}
Since $h_c^\mathrm{que}(\beta) \leq h_c^\mathrm{ann}(\beta)$, this
proves that
$h_c^\mathrm{que}(\beta) = h_c^\mathrm{ann}(\beta)$.
\end{pf}

We now give the proof of Lemma~\ref{supAttained}.

\begin{pf}
The proof is accomplished in three steps. The claims in Steps 1 and~2
are obvious when
the support of $\mu_0$ is bounded from above, because then $\Phi$ is
bounded from above
and upper semicontinuous. Thus, for these steps we may assume that the
support of $\mu_0$
is unbounded from above.\looseness=1

\textit{Step} 1: The supremum can be restricted to the set $\cC\cap\{
Q\in\cP^\mathrm{inv}
(\widetilde E^\N)\dvtx \break I^\mathrm{que}(Q)\leq\gamma\}$ for some
$\gamma
<\infty$.

\begin{pf}
We first prove that
%
\begin{equation}
\label{EntropyDomination}
\lim_{a\to\infty} \sup_{ {Q\in\cC} \atop{\Phi(Q)=a} }
[\beta\Phi(Q)-I^\mathrm{que}(Q)] = -\infty.
\end{equation}
To that end we estimate, for $a\in(0, \infty)$,
%
\begin{eqnarray}
\qquad\mathop{\sup_{Q\in\cC}}_{\Phi(Q)=a}  [\beta\Phi(Q)-I^\mathrm{que}(Q)]
&\leq&\mathop{\sup_{Q\in\cC}}_{\Phi(Q)=a}
[\beta a-h(\pi_{1,1} Q \mid\mu_0)]
\nonumber
\\[-8pt]
\\[-8pt]
\nonumber
&=&\mathop{\sup_{\mu\in\cP(E)}}_{\int_E |x| \,\di\mu(x) <\infty,
\int_E x
\,\di\mu(x) =a}
[\beta a -h(\mu\mid\mu_0)],
\end{eqnarray}
where we use that $I^\mathrm{que}(Q) \geq I^\mathrm{ann}(Q) = H(Q
\mid
Q_0) \geq
h(\pi_{1,1}Q \mid\mu_0)$. The last~sup\-remum is achieved by a measure
$\mu_\lambda$ of the
form $\di\mu_\lambda(x) = M(\lambda)^{-1}\e^{\lambda x} \,\di\mu_0(x)$,
$x\in E$, with
$\lambda$ such that $\int_E x \,\di\mu_\lambda(x)=a$ [recall (\ref
{mubetadef})]. To see
why, first note that such a $\lambda=\lambda(a)$ exists because
$(\lambda\mapsto\int_E x \,\di\mu_\lambda(x))$ is continuous with value
0 at $\lambda=0$ and
$\lim_{\lambda\to\infty}
\int_E x \,\di\mu_\lambda(x) = \mathrm{\sup}[\mathrm{supp}(\mu_0)]=w$,
where $w=\infty$
by assumption. Next note that, for any other measure $\mu$ with $\int_E
x \,\di\mu(x)=a$,
we have
%
\begin{equation}
h(\mu\mid\mu_\lambda) = h(\mu\mid\mu_0) - \lambda a
+ \log M(\lambda) = h(\mu\mid\mu_0)-h(\mu_\lambda\mid\mu_0),
\end{equation}
which shows that $h(\mu\mid\mu_0) \geq h(\mu_\lambda\mid\mu_0)$ with
equality if
and only if $\mu=\mu_\lambda$. Consequently,
%
\begin{eqnarray}
\qquad\mathop{\sup_{\mu\in\cP(E)}}_{\int_E |x| \,\di\mu(x)<\infty,\int
_E x \,\di
\mu(x)=a}
[\beta a-h(\mu\mid\mu_0)]
&=\phantom{:}& \beta\int_E x \,\di\mu_\lambda(x) - h(\mu_\lambda\mid\mu_0)
\nonumber
\\[-8pt]
\\[-8pt]
\nonumber
 &=:&
g(\lambda).
\end{eqnarray}
Clearly, $a\to\infty$ implies $\lambda=\lambda(a)\to\infty$, and
so to prove
(\ref{EntropyDomination}) we must show that $\lim_{\lambda\to\infty}
g(\lambda)=-\infty$.

To achieve the latter, note that a lower bound on $h(\mu_\lambda\mid
\mu_0)$ is obtained
by applying (\ref{LegendreTransIneq}) to $f(x):=\bar\beta(x\vee0)$
for some $\bar\beta
>\beta$. This yields
%
\begin{equation}
g(\lambda) \leq-(\bar\beta-\beta) \int_E x \,\di\mu_\lambda
(x)+\log
[M(\bar\beta)+1].
\end{equation}
The integral in the right-hand side tends to infinity as $\lambda\to
\infty$, and so
(\ref{EntropyDomination}) indeed follows.

Finally, recall the definition of $A(\beta)$ in \eqref{Adef}, which is
finite because of
Lemma~\ref{ACont}. Then, by (\ref{EntropyDomination}), there is an
$a_0<\infty$ such that
%
\begin{equation}
\mathop{\sup_{Q \in\cC}}_{\Phi(Q)=a}
[\beta\Phi(Q)-I^\mathrm{que}(Q)] \leq A(\beta)-1\qquad \forall a \geq a_0,
\end{equation}
and so all $Q\in\cC$ with $\beta\Phi(Q)-I^\mathrm{que}(Q)>A(\beta)-1$
must satisfy $\Phi(Q)<a_0$
and $I^\mathrm{que}(Q) < \beta\Phi(Q)+1-A(\beta) \leq\beta a_0 +
1-A(\beta)=:\gamma$.
Consequently, the supremum can be restricted to the set $\cC\cap\{
Q\in
\cP^\mathrm{inv}
(\widetilde E^\N)\dvtx  I^\mathrm{que}(Q)\leq\gamma\}$.
\end{pf}

\textit{Step} 2: $\Phi$ is upper semicontinuous on $\{Q\in\cP
^\mathrm
{inv}(\widetilde E^\N)\dvtx
I^\mathrm{que}(Q)\leq\gamma\}$ for every $\gamma>0$.\vadjust{\goodbreak}

\begin{pf}
From the definition of $\Phi$ and the inequality $h(\pi_{1, 1} Q \mid
\mu_0) \leq\break  I^\mathrm{que}(Q)
\leq\gamma$, it follows that it is enough to show that the map $\mu
\mapsto\Psi(\mu) := \int_E\eqref{EntropyDensity},
(x \vee0) \,\di\mu(x)$ is upper semicontinuous on $K_\gamma:=\{\mu
\in\cP
(E)\dvtx\break  h(\mu\mid\mu_0)
\leq\gamma\}$. To do so, let $(\mu^M)_{M\in\N}$ be a sequence in
$K_\gamma$ converging to $\mu$
weakly as $M\to\infty$. Then
%
\begin{equation}
\Psi(\mu^M) =\int_E [ (x \vee0) \wedge n] \,\di\mu^M(x)+ \int_E x 1_{\{x>n\}} \,\di\mu^M(x),
\end{equation}
and so
%
\begin{eqnarray}
\label{split}
\limsup_{M\to\infty} \Psi(\mu^M) &\leq&\int_E [(x \vee0) \wedge
n] \,\di
\mu(x)
\nonumber
\\[-8pt]
\\[-8pt]
\nonumber
&&{}+ \sup_{M\in\N} \int_E x 1_{\{x>n\}} \,\di\mu^M(x)\qquad\forall n\in\N.
\end{eqnarray}
By the inequality in (\ref{LegendreTransIneq}), we have
%
\begin{eqnarray}
\lambda\int_E x 1_{\{x>n\}} \,\di\mu^M(x)\leq h(\mu^M \mid\mu_0)
+ \log\int_E \e^{\lambda x 1_{\{x>n\}}} \,\di\mu_0(x)
\nonumber
\\[-8pt]
\\[-8pt]
\eqntext{\forall M,n\in\N, \lambda>0,}
\end{eqnarray}
and so
%
\begin{equation}
\sup_{M\in\N} \int_E x 1_{\{x>n\}} \,\di\mu^M(x) \leq\frac{\gamma
}{\lambda}
+ \frac{1}{\lambda} \log\int_E \e^{\lambda x 1_{\{x>n\}}} \,\di\mu_0(x).
\end{equation}
By (\ref{mgfcond}), the limit as $n\to\infty$ of the right-hand side is
$\gamma/\lambda$.
Since $\lambda>0$ is arbitrary, we conclude that the limit as $n\to
\infty$ of the left-hand
side is zero. Letting $n\to\infty$ in~(\ref{split}) and using monotone
convergence, we
therefore get $\limsup_{M\to\infty} \Psi(\mu^M) \leq\Psi(\mu)$,
as required.
\end{pf}

\textit{Step} 3: Let $\Gamma(Q):=\beta\Phi(Q)-I^\mathrm{que}(Q)$.
Then, by Step 1, we have that for some $\gamma>0$,
%
\begin{equation}
\sup_{Q\in\cC}\Gamma(Q)
= \mathop{\sup_{Q\in\cC}}_{ I^\mathrm{que}(Q) \leq\gamma} \Gamma(Q)
\leq\mathop{\sup_{Q\in\cP^\mathrm{inv}(\widetilde E^\N)}}_{ I^\mathrm
{que}(Q) \leq\gamma}
\Gamma(Q).
\end{equation}
By Theorem~\ref{qLDP}, $I^\mathrm{que}$ is lower semicontinuous. Hence,
by Step 2,
$\beta\Phi-I^\mathrm{que}$ is upper semicontinuous on the compact set
$\{Q\in\cP^\mathrm{inv}(\widetilde E^\N)\dvtx  I^\mathrm
{que}(Q)\leq
\gamma\}$,
achieving its supremum at some $Q^*$. Let $\mu^*:=\pi_{1,1} Q^*$. Then,
by (\ref{mgfcond}),
the inequality in~(\ref{LegendreTransIneq}) gives
%
\begin{equation}
\int_E (x\vee0) \,\di\mu^*(x) \leq\gamma+\log\int_E \e^x \,\di\mu_0(x)
< \infty,
\end{equation}
and, since $\Phi(Q^*)>-\infty$, we also have $\int_E (x \wedge0)
\,\di\mu
^*(x)>-\infty$, so
that $Q^*\in\cC$. Hence
%
\begin{equation}
\sup_{Q\in\cC} \Gamma(Q)
=\mathop{\sup_{Q\in\cP^\mathrm{inv}(\widetilde E^\N)}}_{ I^\mathrm
{que}(Q)\leq\gamma}
\Gamma(Q) = \Gamma(Q^*),
\end{equation}
which completes the proof.
\end{pf}

\section{Reformulation of the criterion for disorder relevance}
\label{S5}

Note that, by (\ref{eqgndefinitionIalg}) and (\ref{truncapproxcont}), for
$\alpha>0$, the
necessary and sufficient condition for relevance, $I^\mathrm
{que}(Q_\beta)>I^\mathrm{ann}
(Q_\beta)$, in Theorem~\ref{CriterionLemma} translates into
%
\begin{equation}
\label{RelCriterion}
\lim_{\tr\to\infty} m_{[Q_\beta]_{{\tr}}}
H\bigl(\Psi_{[Q_\beta]_{\tr}} \mid\mu_0^{\otimes\N_0}\bigr) > 0.
\end{equation}
In Lemma~\ref{SpRelEntrExpression} below, we give two alternative
expressions for the
specific relative entropy appearing in (\ref{RelCriterion}). These
expressions will
be needed in Sections~\ref{S6} and~\ref{S7}.

\textit{I. Asymptotic mean stationarity}.
In what follows we will make use of the notion of \textit{asymptotic
mean stationarity}
(see Gray~\cite{GR90}, Section 1.7). Let $A$ be a topological space and
equip $A^{\N_0}$
with the product topology. A measure $\cP$ on $A^{\N_0}$ is called
\textit{asymptotically
mean stationary} if for every Borel measurable $G\subset A^{\N_0}$,
%
\begin{equation}
\overline{\cP}(G) := \lim_{n\to\infty} \frac{1}{n} \sum
_{k=0}^{n-1} \cP
(\theta^{-k} G)\qquad
\mbox{exists}.
\end{equation}
As in Section~\ref{S2}, $\theta$ denotes the left-shift acting on
$A^{\N_0}$.
If $\cP$ is asymptotically mean stationary, then $\overline{\cP}$ is a
stationary
measure, called the \textit{stationary mean} of $\cP$.

For $Q\in\cP^\mathrm{inv}(\widetilde{E}^\N)$, recall from
Section~\ref
{S2.1} that
$\kappa(Q)\in\cP(E^{\N_0})$ is the probability measure induced by the
concatenation
map $\kappa\dvtx \widetilde{E}^\N\to E^{\N_0}$ that glues a sequence
of words into
a sequence of letters, that is, $\kappa(Q)= Q \circ\kappa^{-1}$. Our
aim is
to replace $\Psi_Q$ in (\ref{RelCriterion}) by $\kappa(Q)$, which is
not stationary but
more convenient to work with. These two probability measures are
related in the
following way.

\begin{itlemma} \label{AMSLemma}
If $m_Q<\infty$, then $\kappa(Q)$ is asymptotically mean stationary
with stationary mean
$\overline{\kappa(Q)}=\Psi_Q.$
\end{itlemma}

\begin{pf}
Let $X:=\kappa(Y) \in E^{\N_0}$, where $Y$ is distributed according to
$Q$. Let $I$ denote
the set of indices $i\in\N_0$ where a new word starts ($0\in I$). For
$i\in\N_0$, let
$r_i:=\inf\{j\in\N\dvtx  i-j\in I\}$, that is, the distance from $i$ to
the beginning of
the word it belongs to. For $j \in I$, let $L^j$ denote\vadjust{\goodbreak} the length of
the word that
starts at $j$. Then, for any $G\subset E^{\N_0}$ Borel measurable, we
have
%
\begin{eqnarray}
\label{sum11}
\sum_{i=0}^{n-1} \kappa(Q)(\theta^i X \in G)
&=&\sum_{i=0}^{n-1} \sum_{k=0}^i Q(\theta^i X\in G, r_i=k)
\nonumber
\\[-8pt]
\\[-8pt]
\nonumber
&=& \sum_{k=0}^{n-1} \sum_{i=k}^{n-1} Q(\theta^i X \in G,
r_i=k).
\end{eqnarray}
Next, note that
%
\begin{eqnarray}
&&Q(\theta^i X \in G, r_i=k)\nonumber\\
&&\qquad= Q(\theta^i X \in G, i-k\in I, L^{i-k}>k)
\nonumber
\\[-8pt]
\\[-8pt]
\nonumber
&&\qquad= Q(\theta^i X \in G, L^{i-k}>k \mid i-k\in I) Q(i-k\in I)\\
&&\qquad= Q(\theta^k X \in G, L^0>k) Q(i-k\in I).\nonumber
\end{eqnarray}
Hence, dividing the sum in (\ref{sum11}) by $n$, we get
%
\begin{equation}
\frac{1}{n} \sum_{i=0}^{n-1} \kappa(Q)(\theta^i X \in G)
= \sum_{k=0}^{n-1} Q(\theta^k X \in G, L^0>k) f_{k,n},
\end{equation}
where we abbreviate $f_{k,n}:=n^{-1}\sum_{j=0}^{n-k-1} Q(j\in I)$. By
the renewal theorem,
$\lim_{n\to\infty} f_{k, n}= 1/m_Q$ for $k$ fixed. Since
%
\begin{equation}
\sum_{k=0}^\infty Q(L^0>k)=m_Q<\infty,
\end{equation}
we can apply the bounded convergence theorem, and conclude that
%
\begin{eqnarray}
\overline{\kappa(Q)}(G)
&=&\frac{1}{m_Q} \sum_{k=0}^\infty Q(\theta^k X \in G, L^0>k)
\nonumber\\
&=&\frac{1}{m_Q} \sum_{k=0}^\infty\sum_{j=k+1}^\infty Q(\theta^k X
\in
G, L^0=j)
\\
&=&\frac{1}{m_Q} \sum_{j=1}^\infty\sum_{k=0}^{j-1} Q(\theta^k X \in
G, L^0=j)
=\Psi_Q(G).\nonumber
\end{eqnarray}
The last equality is simply the definition of $\Psi_Q$ in (\ref{PsiQdef}).
\end{pf}

To complement Lemma~\ref{AMSLemma}, we need the following fact stated
in Birkner~\cite{Bi08},
Remark 5, where ergodicity refers to the left-shifts acting on
$\widetilde{E}^\N$ and $E^\N$.

\begin{itlemma} \label{BirknerRemark}
If $Q\in\cP^\mathrm{inv}(\widetilde{E}^\N)$ is ergodic and
$m_Q<\infty
$, then $\Psi_Q
\in\break\cP^\mathrm{inv}(E^\N)$ is ergodic.\vadjust{\goodbreak}
\end{itlemma}

An asymptotic mean stationary measure can be interchanged with its
stationary mean
in several situations (see Gray~\cite{GR88}, Chapter 6), for example, in
relative entropy computations, as in Lemma~\ref{SpRelEntrExpression}
below. Before stating
this lemma, we use an extension of the notion of specific relative
entropy to measures that
are not necessarily stationary. More precisely, for two measures $\cP$
and $\cQ$ on a product
space~$A^\N$, we define the specific relative entropy of $\cP$ w.r.t.
$\cQ$ as
%
\begin{equation}
\overline{H}(\cP\mid\cQ) := \limsup_{n\to\infty} \frac{1}{n}
h(\pi_n \cP\mid\pi_n \cQ),
\end{equation}
where $\pi_n$ is the projection onto the first $n$ coordinates. For
$Q\in\cP^\mathrm{inv}
(\widetilde{E}^\N)$, we introduce the following Radon--Nikodym derivative:
%
\begin{equation}
\label{fnRN}
f_n(x):=\frac{\di\pi_n\kappa(Q)}{\di\mu_0^{\otimes n}} (x),\qquad
x\in E^{\N_0}.
\end{equation}
With this notation, the main result of this section is the following.

\begin{itlemma}
\label{SpRelEntrExpression}
For $Q\in\cP^\mathrm{inv}(\widetilde{E}^\N)$ ergodic with
$m_Q<\infty$,
\begin{eqnarray}
\label{RelEntr*}\quad
H(\Psi_Q\mid\mu_0^{\otimes\N_0})
&=& \overline{H}(\kappa(Q)\mid\mu_0^{\otimes\N_0}),\\
\label{RelEntr}
&= &\lim_{n\to\infty} \frac{1}{n}
\log f_n(x)\qquad \mbox{for $\kappa(Q)$-a.s. all $x\in E^{\N_0}$}.
\end{eqnarray}
The first equality holds also without the assumption of ergodicity.
\end{itlemma}

\begin{pf}
The first equality follows from Gray~\cite{GR90}, Corollary 7.5.1, last
equality in
equation (7.32), which does not need the assumption of ergodicity. For
the proof of
the other equality, define
%
\begin{equation}
\bar{f}_n(x):=\frac{\di\pi_n\Psi_Q}{\di\mu_0^{\otimes n}} (x).
\end{equation}
Since $\Psi_Q$ is stationary and ergodic (Lemma~\ref{BirknerRemark}),
Gray~\cite{GR90},
Theorem 8.2.1, applied to the pair $\Psi_Q$, $\mu_0^{\otimes\N_0}$
gives that
%
\begin{equation}
\label{RNConvergence}
\lim_{n\to\infty}\frac{1}{n}\log\bar{f}_n(x)=H(\Psi_Q\mid\mu
_0^{\otimes
\N_0})
\end{equation}
for $\Psi_Q$ almost all $x$. But $\Psi_Q$ is the stationary mean of
$\kappa(Q)$
(Lemma~\ref{AMSLemma}), so that Gray~\cite{GR90}, Theorem 8.4.1,
combined with
\eqref{RNConvergence} gives
%
\begin{equation}
\lim_{n\to\infty}\frac{1}{n}\log f_n(x)=H(\Psi_Q\mid\mu
_0^{\otimes\N_0})
\end{equation}
for $\kappa(Q)$ almost all $x$.
\end{pf}

\textit{II. Alternative formulation}.
We will apply Lemma~\ref{SpRelEntrExpression} to the measure $[Q_\beta
]_{\tr}$, which is
ergodic, being a product measure.\vadjust{\goodbreak} The word length distribution of it is
%
\begin{equation}
\label{Ktrun}
K^{\tr}(n):=
\cases{ K(n), & \quad $\mbox{if } 1 \leq n \leq\tr-1,$\vspace*{2pt}\cr
\displaystyle\sum_{m=\tr}^\infty K(m), &\quad $\mbox{if } n=\tr,$\vspace*{2pt}\cr
0 ,& \quad$\mbox{if } n>\tr.$}
\end{equation}
For $[Q_\beta]_{\tr}$, the function $f_n$ in (\ref{fnRN}) becomes
%
\begin{equation}
\label{EntropyDensity}
\qquad f_n(x)=\ev_{K^{\tr}}\Biggl(\prod_{k=0}^{n-1}
\biggl(\frac{\e^{\beta x_k}}{M(\beta)}\biggr)^{1_{\{S_k=0\}}}\Biggr)
=\ev_{K^{\tr}}\bigl(\e^{\sum_{k=0}^{n-1}
\{\beta x_k-\log M(\beta)\}1_{\{S_k=0\}}}\bigr),\vspace*{6pt}
\end{equation}
where $\ev_{K^{\tr}}$ denotes expectation with respect to law of the
Markov chain $S$
with renewal time distribution $K^{\tr}$ starting from 0. This follows
from the definition
of $Q_\beta$ and (\ref{mubetadef}). To emphasize the fact that in the
last expression
the sequence $x\in E^{\N_0}$ is picked from $\kappa([Q_\beta]_{\tr}
)$, we
take two independent
sequences
%
\begin{equation}
\label{xtilt}\quad
(x_k)_{k\in\N_0}, (\hat x_k)_{k\in\N_0} \mbox{ drawn from $\mu
_0^{\otimes\N_0}$
and $\mu_\beta^{\otimes\N_0}$, respectively},
\end{equation}
and an independent copy $S'$ of $S$. Let $I:=\{i\ge0: S_i=0\}, I':=\{
i\ge0:\break S'_i=0\}$. Then
%
\begin{eqnarray}
\label{SRelativeEntropyAlternative}
&& H\bigl(\Psi_{[Q_\beta]_{\tr}}
\mid\mu_0^{\otimes\N_0}\bigr)
\nonumber
\\[-8pt]
\\[-8pt]
\nonumber
&&\qquad =\lim_{n\to\infty}\frac{1}{n}\log\ev
_{K^{\tr}}
\bigl[\e^{ \sum_{k=0}^{n-1}[\beta x_k 1_{\{k\notin I'\}}+\beta\hat x_k
1_{\{k\in I'\}}-\log M(\beta)] 1_{\{k\in I\}}}\bigr].
\end{eqnarray}
Note the appearance of two renewal sets $I,I'$, which are the key to
understanding the
issue of relevant vs. irrelevant disorder; recall Remark~\ref{I1I2}.


\section{\texorpdfstring{Monotonicity of disorder relevance: Proof of Theorem~\protect\ref{MonRelevance}}
{Monotonicity of disorder relevance: Proof of Theorem 1.6}}
\label{S6}
\mbox{}
\begin{pf}
In view of (\ref{RelEntr*}) in Lemma~\ref{SpRelEntrExpression}, the
condition for
relevance in (\ref{RelCriterion}) becomes
%
\begin{equation}
\label{RelCriterion2}
\lim_{\tr\to\infty} m_{[Q_\beta]_{\tr}}
\overline{H}(\kappa([Q_\beta]_{\tr}) \mid\mu_0^{\otimes\N_0}) > 0.
\end{equation}
We will show that $\beta\mapsto\overline{H}(\kappa([Q_\beta]_{\tr})
\mid\mu_0^{\otimes\N_0})$
is nondecreasing for every $\tr\in\N$, which will imply the claim because
$m_{[Q_\beta]_{\tr}}=m_{K^{\tr}}$ does not depend on $\beta$. It will be
enough to show that
$\beta\mapsto h(\pi_n \kappa([Q_\beta]_{\tr}) \mid\mu_0^{\otimes n})$
is nondecreasing
for all $\tr,n\in\N$.

Fix $\tr,n\in\N$. For $\beta\in[0,\infty)$ and $\bar{x}=(x_0, x_1,
\ldots, x_{n-1})\in E^n$, let
%
\begin{equation}
k(\beta,\bar{x}) := \frac{\di\pi_n \kappa([Q_\beta]_{\tr})}{\di
\mu
_0^n}(\bar{x})
= \ev_{K^{\tr}} \biggl(\prod_{k\in J_n} \frac{\e^{\beta x_k}}{M(\beta)}\biggr),
\end{equation}
with $J_n :=\{0\leq k < n\dvtx  S_k=0\}$ the set of renewal times
prior to time $n$ for the
chain $S$ that has renewal time distribution $K^{\tr}$, to which we add
$0$ for convenience. Our
goal is to prove that
%
\begin{equation}
\qquad\beta\mapsto
f(\beta) := \int_{\R^n} [k(\beta,\bar{x}) \log k(\beta,\bar{x})]
\,\di\mu^{\otimes n}_0(\bar{x}) = h(\pi_n\kappa([Q_\beta]_{\tr})
\mid\mu
_0^{\otimes n})
\end{equation}
is nondecreasing on $[0,\infty)$. We will do this by proving a stronger
property. Namely,
for $\bar{\beta}=(\beta_0, \beta_1, \ldots, \beta_{n-1}) \in
[0,\infty
)^n$ and $\bar{x}\in E^n$,
let
%
\begin{equation}
k(\bar{\beta},\bar{x}) := \ev_{K^{\tr}} \biggl(\prod_{k\in J_n}
\frac{\e^{\beta_k x_k}}{M(\beta_k)}\biggr).
\end{equation}
We will show that
%
\begin{equation}
\bar{\beta} \mapsto f(\bar{\beta})
:= \int_{\R^n} [k(\bar{\beta},\bar{x})
\log k(\bar{\beta},\bar{x})] \,\di\mu^{\otimes n}_0(\bar{x})
\end{equation}
is nondecreasing on $[0,\infty)^n$ in each of its arguments.

We will prove monotonicity w.r.t. $\beta_1$ only. The argument is the
same for the other
variables, with one simplification for $\beta_0$; namely, we may drop
the corresponding
indicator $1_{\{0\in J_n\}}$ in the third line of (\ref
{entropyDerivative}) and in
(\ref{x1Function}). First, using that $\int k(\bar{\beta},\bar
{x})\,\di\mu
_0^{\otimes n}
(\bar{x})=1$ for all $\bar{\beta}$, we compute
%
\begin{eqnarray}
\label{entropyDerivative}
&&\partial_{\beta_1} f(\bar{\beta})\nonumber\hspace*{-35pt}\\
&&\quad= \int_{\R^n} \partial_{\beta_1} [k(\bar{\beta},\bar{x})
\log k(\bar{\beta},\bar{x})] \,\di\mu^{\otimes n}_0(\bar{x})
\nonumber\hspace*{-35pt}
\\[-8pt]
\\[-8pt]
\nonumber
&&\quad= \int_{\R^n} \partial_{\beta_1} [k(\bar{\beta},\bar{x})]
\log k(\bar{\beta},\bar{x}) \,\di\mu^{\otimes n}_0(\bar{x})\hspace*{-35pt}\\
&&\quad= \int_{\R^n} \partial_{\beta_1} \biggl(\frac{\e^{\beta_1
x_1}}{M(\beta_1)}\biggr)
\ev_{K^{\tr}} \biggl(1_{\{1\in J_n\}} \prod_{k\in J_n\setminus\{1\}}
\frac{\e^{\beta_k x_k}}{M(\beta_k)} \biggr)
\log k(\bar{\beta},\bar{x}) \,\di\mu^{\otimes n}_0(\bar{x}).\nonumber\hspace*{-35pt}
\end{eqnarray}
Next, we note that
%
\begin{eqnarray}
\partial_{\beta_1} \biggl(\frac{\e^{\beta_1 x_1}}{M(\beta_1)}\biggr) \,\di\mu_0(x_1)
&=& \frac{\e^{\beta_1 x_1} x_1 M(\beta_1) -
\e^{\beta_1 x_1}M'(\beta_1)}{M(\beta_1)^2} \,\di\mu_0(x_1)
\nonumber
\\
&=& \biggl(x_1-\frac{M'(\beta_1)}{M(\beta_1)}\biggr)
\frac{\e^{\beta_1 x_1}}{M(\beta_1)} \,\di\mu_0(x_1)\\
&=&(x_1-E_{\beta_1}) \,\di\mu_{\beta_1}(x_1),\nonumber
\end{eqnarray}
where $E_{\beta_1}:=M'(\beta_1)/M(\beta_1) =\int x_1 \,\di\mu_{\beta
_1}(x_1)$.
Now, let $\bar{x}^1$ be $\bar{x}$ without $x_1$, and abbreviate
%
\begin{equation} \label{x1Function}
A(x_1;\bar{x}^1) := \ev_{K^{\tr}} \biggl(\prod_{k\in J_n\setminus\{1\}}
\frac{\e^{\beta_k x_k}}{M(\beta_k)} 1_{\{1\in J_n\}}\biggr)
\log k(\bar{\beta},\bar{x}).
\end{equation}
Then, for fixed $\bar{x}^1$, the integral over $x_1$ in (\ref
{entropyDerivative}) equals
%
\begin{eqnarray}
&&\int_{\R^n} (x_1-E_{\beta_1}) A(x_1;\bar{x}^1) \,\di\mu_{\beta
_1}(x_1)
\nonumber
\\[-8pt]
\\[-8pt]
\nonumber
&&\qquad \geq\int_{\R^n} (x_1-E_{\beta_1}) \,\di\mu_{\beta_1}(x_1)
\int_{\R^n} A(x_1;\bar{x}^1) \,\di\mu_{\beta_1}(x_1) = 0,
\end{eqnarray}
where the inequality holds because both $x_1\mapsto x_1-E_{\beta_1}$
and $x_1 \mapsto
A(x_1;\bar{x}^1)$ are nondecreasing [for the latter we need that
$\beta
_1\in[0,\infty)$].
It therefore follows from~(\ref{entropyDerivative}), after integrating
over $\bar{x}^1$
as well, that $\partial_{\beta_1} f(\bar{\beta})\geq0$.
\end{pf}


\section{\texorpdfstring{Disorder irrelevance: Proof of Corollaries~\protect\ref{disIrrelevancealphazero} and \protect\ref{BdsCritTemp}(i)}
{Disorder irrelevance: Proof of Corollaries 1.7 and 1.8(i)}}
\label{S7}


\subsection{\texorpdfstring{Proof of Corollary~\protect\ref{disIrrelevancealphazero}}
{Proof of Corollary 1.7}}
\label{S7.1}
\mbox{}
\begin{pf}
This is immediate from Theorem~\ref{CriterionLemma} and the fact that
$I^\mathrm{que}
=I^\mathrm{ann}$ when $\alpha=0$. The latter was already noted at the
end of Section~\ref{S2}.
\end{pf}


\subsection{\texorpdfstring{Proof of Corollary~\protect\ref{BdsCritTemp}\textup{(i)}}
{Proof of Corollary 1.8(i)}}
\label{S7.2}
\mbox{}
\begin{pf}We will show disorder irrelevance for all $\beta$ that satisfy
$M(2\beta)/\break M(\beta)^2<1
+\chi^{-1}$. To show that for such $\beta$ the limit in \eqref
{RelCriterion} is zero, we
use an annealed bound on $H(\Psi_{[Q_\beta]_{\tr}} \mid\mu
_0^{\otimes\N
_0})$ based on the
expression \eqref{RelEntr} for it. We bound the limit in the right-hand
side of that formula,
using \eqref{annealedBound} with the role of $\Theta_n$ played by
%
\begin{equation}
\label{fnxdef}
f_n(x)=\frac{\di\pi_n \kappa([Q_\beta]_{\tr})}{\di\mu_0^{\otimes n}}(x),
\qquad x\in E^{\N_0}.
\end{equation}
This satisfies
%
\begin{equation}
\label{fnxrel}
\ev_{\kappa([Q_\beta]_{\tr})}(f_n(x)) = \ev_{\mu_0^{\otimes
n}}(f_n(x) f_n(x)),
\end{equation}
because $f_n(x)$ depends on the first $n$ coordinates of $x$ only, and
the Radon--Nikodym
derivative of $\pi_n\kappa([Q_\beta]_{\tr})$ with respect to $\mu
_0^{\otimes n}$ is $f_n$.
Using \eqref{EntropyDensity}, we write the last expectation as
%
\begin{eqnarray}
\label{replica}
&&\ev_{\mu_0^{\otimes n}}(f_n(x) f_n(x))\nonumber\\
&&\qquad= \ev_{\mu_0^{\otimes n}}\Biggl((\ev_{K^{\tr}}\times\ev_{K^{\tr}})
\Biggl(
\prod_{k=0}^{n-1}\biggl(\frac{\e^{\beta x_k}}{M(\beta)}\biggr)^{1_{\{S_k=0\}}}
\prod_{l=0}^{n-1}\biggl(\frac{\e^{\beta x_l}}{M(\beta)}\biggr)^{1_{\{S'_l=0\}}}
\Biggr)\Biggr)
\nonumber
\\[-8pt]
\\[-8pt]
\nonumber
&&\qquad=(\ev_{K^{\tr}}\times\ev_{K^{\tr}})
\Biggl(\ev_{\mu_0^{\otimes n}}\Biggl(
\prod_{k=0}^{n-1}\biggl(\frac{\e^{\beta x_k}}{M(\beta)}\biggr)^{1_{\{S_k=0\}}}
\prod_{l=0}^{n-1}\biggl(\frac{\e^{\beta x_l}}{M(\beta)}\biggr)^{1_{\{S'_l=0\}}}
\Biggr)\Biggr)\\
&&\qquad=(\ev_{K^{\tr}}\times\ev_{K^{\tr}})\bigl(\Xi(\beta)^{\sum
_{k=0}^{n-1} 1_{\{
S_k=S'_k=0\}}}\bigr),\nonumber
\end{eqnarray}
where $\ev_{K^{\tr}}\times\ev_{K^{\tr}}$ is the expectation with
respect to two independent
copies $S,S'$ of the Markov chain starting from 0 with renewal time
distribution $K^{\tr}$, and
%
\begin{equation}
\label{Xibetadef}
\Xi(\beta) :=\frac{M(2\beta)}{M(\beta)^2}.
\end{equation}
If we now let
%
\begin{equation} \label{trFreeEnergy}
f_2^{\tr}(\lambda) := \lim_{n\to\infty} \frac{1}{n}
\log(\ev_{K^{\tr}}\times\ev_{K^{\tr}})
\bigl(\e^{\lambda\sum_{k=0}^{n-1} 1_{\{S_k=S_k'=0\}}}\bigr),
\end{equation}
then \eqref{RelEntr}, \eqref{annealedBound} and (\ref{fnxdef})--(\ref
{trFreeEnergy}) imply that
%
\begin{equation}
\label{SREUpperBound}
H\bigl(\Psi_{[Q_\beta]_{\tr}} \mid\mu_0^{\otimes\N_0}\bigr) \leq f_2^{\tr}
(\log\Xi
(\beta)),\qquad
\beta\in[0,\infty), \tr\in\N.
\end{equation}
Combining this bound with the condition for relevance in (\ref
{RelCriterion}), we see that to prove irrelevance it suffices to show that
%
\begin{equation}
\label{zeroLimit}
\lim_{\tr\to\infty} m_{[Q_\beta]_{\tr}} f_2^{\tr}(\log\Xi(\beta
)) = 0.
\end{equation}
By (\ref{fhomozerocond}) in \hyperref[A]{Appendix~A}, we have
%
\begin{equation}
\label{zeroLimithomo}
f_2(\lambda) = 0 \quad\Longleftrightarrow\quad
\lambda\leq\lambda_0 := -\log\pr(I\cap I'\neq\varnothing),
\end{equation}
where $I,I'$ are the sets of renewal times for $S,S'$ without
truncation, and~$f_2(\lambda)$ as defined in \hyperref[A]{Appendix~A}. By Lemma~\ref
{freeEnergyTruncated},
if $\lambda<\lambda_0$, then\break $\sup_{tr\in\N}\tr f_2^{\tr}(\lambda)$
$<\infty$. Since
$\lim_{\tr\to\infty} m_{[Q_\beta]_{\tr}}/\tr=0$ always, (\ref
{zeroLimit}) holds as
soon as $\log\Xi(\beta)<\lambda_0$, that is, $\Xi(\beta)<1/\pr
(I\cap
I'\neq\varnothing)$.
Now the claim of the corollary follows because $\pr(I\cap I'\neq
\varnothing)=\chi/(\chi+1)$ (see
Spitzer~\cite{Sp76}, Section~1), with $\chi$ as defined in~\eqref
{chiDefinition},
and with the convention that the last ratio is 1 if $\chi=\infty$.
\end{pf}


\section{\texorpdfstring{Disorder relevance: Proof of Corollary~\protect\ref{BdsCritTemp}(ii)}
{Disorder relevance: Proof of Corollary 1.8(ii)}}
\label{S8}
\mbox{}
\begin{pf}
We restrict the expectation in \eqref{SRelativeEntropyAlternative} to
the set
%
\begin{equation}
A_n:=\bigl\{(S_k)_{k=0}^n\dvtx  I\cap\{1,\ldots, n\}=I'\cap\{1,\ldots,n\}
\bigr\},
\end{equation}
that is, $S$ follows $I'$ and collects only the tilted charges $\hat
x_k$ defined in
(\ref{xtilt}). This gives for the expectation the lower bound
%
\begin{equation}
\exp\Biggl[\sum_{k=0}^{n-1} [\beta\hat x_k-\log M(\beta)] 1_{\{k\in I'\}}\Biggr]
\pr(A_n).
\end{equation}
Let $k_n:=|I\cap\{1,\ldots, n\}|$, $\tau_0'=0$ and $\tau_1'<\cdots
<\tau
_{k_n}'$ the
elements of $I'\cap\{1,\ldots,n\}$. By the renewal theorem, we have
$k_n/n\to1/m_{\tr}$
as $n\to\infty$. Moreover,
%
\begin{equation}
\pr(A_n) = \pr(\tau_1>n-\tau_{k_n}')\prod_{i=1}^{k_n} K^{\tr}(\tau
_i'-\tau
_{i-1}'),
\end{equation}
so that
%
\begin{eqnarray}
\frac{1}{n}\log\pr(A_n)&=&\frac{1}{n} \log\pr(\tau_1>n-\tau
_{k_n})+\frac{k_n}{n}
\frac{1}{k_n}\sum_{i=1}^{k_n}\log K^{\tr}(\tau_i'-\tau_{i-1}')
\nonumber
\\[-8pt]
\\[-8pt]
\nonumber
&\to&
\frac
{1}{m_{\tr}}
\sum_{k=1}^{\tr} K^{\tr}(k)\log K^{\tr}(k),
\end{eqnarray}
while
%
\begin{equation}
\frac{1}{n}\sum_{k=0}^{n-1}\{\beta\hat x_k-\log M(\beta)\} 1_{\{
k\in
I'\}}
\to\frac{1}{m_{\tr}} c(\beta)
\end{equation}
with
%
\begin{equation}
\qquad c(\beta):=\beta E_{\mu_\beta}(\hat x_1)-\log M(\beta)=\beta[\log
M(\beta)]'
-\log M(\beta)=h(\mu_\beta\mid\mu_0).
\end{equation}
Hence
%
\begin{equation}
m_{\tr} H\bigl(\Psi_{[Q_\beta]_{\tr}}\mid\mu_0^{\otimes\N_0}\bigr)
\geq h(\mu_\beta\mid\mu_0) +\sum_{k=0}^{\tr} K^{\tr}(k)\log K^{\tr}(k),
\end{equation}
and
%
\begin{equation}
\liminf_{\tr\to\infty}
m_{[Q_\beta]_{\tr}} H(\kappa([Q_\beta]_{\tr}) \mid\mu_0^{\otimes\N_0})
\geq h(\mu_\beta\mid\mu_0)-H(K).
\end{equation}
Consequently, $h(\mu_\beta\mid\mu_0)>H(K)$ is sufficient for
disorder relevance.
\end{pf}

We close by proving the second part of (\ref{twolims}).
%
\begin{equation}
\label{twolimsalt}
\lim_{\beta\to\infty} h(\mu_\beta\mid\mu_0) = \log[1/\mu_0(\{
w\})].
\end{equation}
We distinguish three different cases:

(1) $w=\infty$. Apply (\ref{LegendreTransIneq}) with $\mu=\mu_\beta$,
$\nu=\mu_0$ and
$f(x)= x \vee0$, to get
%
\begin{equation}
h(\mu_\beta\mid\mu_0) \geq\int_E (x\vee0) \,\di\mu_\beta(x) -
\log[M(1)+1].
\end{equation}
The integral diverges as $\beta\to\infty$, and so (\ref
{twolimsalt}) follows.

(2) $\mu_0(\{w\})=0$ with $w<\infty$. Now $\mu_\beta$ converges weakly
as $\beta\to\infty$
to $\delta_w$, the point measure at $w$. Hence (\ref{twolimsalt})
follows by using the
lower semicontinuity of $\mu\mapsto h(\mu\mid\mu_0)$ and the fact that
$h(\delta_w\mid\mu_0)
=\infty$ because $\delta_w$ is not absolutely continuous w.r.t. $\mu_0$.

(3) $\mu_0(\{w\})>0$ with $w<\infty$. Define
%
\begin{equation}
f_\beta(x) := \frac{\di\mu_\beta}{\di\mu_0}(x)=\frac{\e^{\beta
x}}{M(\beta)},\qquad
x\in E.
\end{equation}
This function satisfies
%
\begin{eqnarray}
\lim_{\beta\to\infty} f_\beta(x)&=&0\qquad \mbox{for } x<w,\nonumber\\
\lim_{\beta\to\infty} f_\beta(w)& =&1/\mu_0(\{w\}),\\
f_\beta(x) &\leq& 1/\mu_0(\{w\})<\infty\qquad \mbox{for } x \leq
w.\nonumber
\end{eqnarray}
Since $t\mapsto t\log t$ is increasing on $[1,\infty)$ and on $(0,1]$
takes values
in $[-\e^{-1},0]$, we can apply the bounded convergence theorem to the integral
%
\begin{equation}
h(\mu_\beta\mid\mu_0) = \int_E f_\beta(x)\log f_\beta(x) \,\di\mu_0(x),
\end{equation}
to get (\ref{twolimsalt}).


\begin{appendix}

\section{Standard facts about the homopolymer}
\label{A}

In this appendix we recall a few standard facts about the homopolymer.
For proofs
we refer to Giacomin~\cite{Gi07}, Chapter 2, and den Hollander~\cite
{dHo09}, Chapter 7.

The homopolymer has a path measure as in (\ref{Pndef}), but with
exponent $\lambda
\sum_{k=0}^{n-1} 1_{\{S_k=0\}}$, $\lambda\in[0,\infty)$. For a given
renewal time
distribution $K$, it is known that the free energy $f(\lambda)$ is the
unique solution
of the equation
%
\begin{equation}
\label{fhomorel}
\mathrm{e}^{-\lambda}=\sum_{n\in\N} K(n) \e^{-nf(\lambda)}
\end{equation}
whenever a solution exists, otherwise $f(\lambda)=0$. Clearly
%
\begin{equation}
\label{fhomozerocond}
f(\lambda)=0 \quad\Longleftrightarrow\quad\lambda\leq-\log\pr(I\neq
\varnothing),
\end{equation}
where $I=\{k\in\N\dvtx  S_k=0\}$ is the set of renewal times of $S$.

Let $S,S'$ be two independent copies of the Markov chain starting form
0, with renewal
time distribution $K$, and with sets of renewal times $I,I'$.
Transience of the joint
renewal process $I\cap I'$ is equivalent to \mbox{$\pr(I\cap I'\neq
\varnothing
)<1$}. In that case,
let
%
\begin{equation}
\lambda_0 := -\log\pr(I\cap I'\neq\varnothing)>0,
\end{equation}
and denote by $f_2(\lambda)$ and $f_2^{\tr}(\lambda)$ the free energy
of the homopolymer
whose underlying Markov chain has renewal set $I \cap I'$ when the
renewal times of $S,S'$
are drawn from $K$, respectively, $K^{\tr}$ defined in (\ref{Ktrun}).
Then\break $\lim_{\tr\to\infty}
f_2^{\tr}(\lambda)=f_2(\lambda)$. Note that $f_2(\lambda)=0$ if and
only if $\lambda\leq\lambda_0$.
This property does not hold for $f_2^{\tr}(\lambda)$, but the following
lemma shows that
$f_2^{\tr}(\lambda)$ tends to zero fast as $\tr\to\infty$ when
$\lambda
<\lambda_0$.

\begin{itlemma}
\label{freeEnergyTruncated}
Suppose that $\pr(I \cap I' \neq\varnothing)<1$. Then $\sup
_{\tr
\in\N}\tr
f_2^{\tr}(\lambda)<\infty$ for all $\lambda<\lambda_0$.
\end{itlemma}

\begin{pf}
As in the paragraph preceding the lemma, define $I^{\tr}, I'^{\tr}$,
where now the Markov
chains $S, S'$ have renewal time distribution $K^{\tr}$. Let $K_2,
K_2^{\tr}
$ be the renewal
time distributions generating the sets $I\cap I', I^{\tr}\cap I'^{\tr}$
respectively. Put
$L_2(n):=\sum_{k=1}^n K_2(k)$ and $L_2^{\tr}(n):=\sum_{k=1}^n K_2^{\tr}
(k)$. Then $L_2(\infty)
=\e^{-\lambda_0}$ and $L_2^{\tr}(\infty)=1$ because the renewal process
$I^{\tr}\cap I'^{\tr}$
is resurrent. Since $K_2^{\tr}(n)=K_2(n)$ for $1\leq n<{\tr}$, it follows
from (\ref{fhomorel})
that
%
\begin{eqnarray}
e^{-\lambda} &=& \sum_{n=1}^{\tr-1} K_2(n) \e^{-n f_2^{\tr}(\lambda)}
+ \sum_{n=\tr}^\infty K_2^{\tr}(n) \e^{-n f_2^{\tr}(\lambda)}
\nonumber
\\[-8pt]
\\[-8pt]
\nonumber
&\leq &L_2(\tr-1) + \e^{-\tr f_2^{\tr}(\lambda)}[1-L_2(\tr-1)],
\end{eqnarray}
where the equality holds because $f_2^{\tr}(\lambda)>0$ for $\lambda
>0$. Hence
%
\begin{equation}
\tr f_2^{\tr}(\lambda) \leq\log\biggl[
\frac{1-L_2(\tr-1)}{\e^{-\lambda}-L_2(\tr-1)}\biggr].
\end{equation}
The term between brackets tends to $(1-\e^{-\lambda_0})/(\e
^{-\lambda
}-\e^{-\lambda_0})$ as
$\tr\to\infty$, which is finite for $\lambda<\lambda_0$.
\end{pf}

The order of the phase transition for the homopolymer depends on the
tail of $K$. If $K$
satisfies (\ref{Kcond}), then (see \cite{Gi07}, Theorem 2.1, \cite
{dHo09}, Theorem 7.4)
%
\begin{equation}
f(\lambda) \sim\lambda^{1/(1 \wedge\alpha)} L^*(1/\lambda),\qquad
\lambda\downarrow0,
\end{equation}
for some $L^*$, that is, strictly positive and slowly varying at
infinity. Hence, the
phase transition is order $1$ when $\alpha\in[1,\infty)$ and order
$m\in\N\setminus
\{1\}$ when $\alpha\in[\frac{1}{m},\frac{1}{m-1})$. This shows that
the value
$\alpha=\tfrac12$ is critical in view of the Harris criterion
mentioned in
Remark~\ref{RHarris}.


\section{\texorpdfstring{Proof of Lemma~\lowercase{\protect\ref{lsc_lemma}}}
{Proof of Lemma 3.2}}
\label{B}

We borrow ideas from the proof of the lower bound of the LDP in
Theorem~\ref{qLDP}
given in Birkner, Greven and den Hollander~\cite{BiGrdHo09},
Proposition 4.1. What
follows is a rewriting of the relevant parts of that proof, organized as
Sections~\ref{B1}--\ref{B4}. Our setting is the same as their setting
because the assumption throughout Section~\ref{S3.2} is that
$E:=\mathrm
{supp}[\mu_0]$ is finite.

We will prove that $S^\mathrm{que}(\beta;1)\le S^\mathrm{que}(\beta
;1-)$. Fix $A<S^\mathrm{que}(\beta;1)$. By \eqref{Sformulaz=1} and~\eqref{truncapproxcont}, there
is a $Q\in\cP^\mathrm{inv}(\widetilde{E}^\N)$ with $m_{Q}<\infty$ such
that $\beta\Phi(Q)
-  I^\mathrm{que}(Q)>A$. Because $\Phi$ and $I^\mathrm{que}$ are affine,
we may assume without
loss of generality that~$Q$ is ergodic.


\subsection{Step 1: Good sentences}
\label{B1}

For $\varepsilon>0$, the set
%
\begin{equation}
\label{QNeighborhood}
\cU_\varepsilon(Q) := \{Q'\in\cP^\mathrm{inv}(\widetilde{E}^\N
)\dvtx
\Phi(Q')>\Phi(Q)-\varepsilon\}
\end{equation}
is open because $\Phi$ is continuous. Hence there is an $M_0\in\N$
large enough, a~$\delta_1>0$
and a finite set $\mathcal{A}_0\subset\widetilde E^{M_0}$ such that
%
\begin{equation}
\mathcal{U}_{M, \delta_1} := \{Q'\in\cP^{\mathrm{inv}}(\widetilde
{E}^\N
)\dvtx
|(\pi_{M_0}Q')(s)-f_s|<2\delta_1\ \forall s\in\mathcal{A}_0\}
\subset\cU_\varepsilon(Q),\hspace*{-35pt}
\end{equation}
where we set $f_s:=(\pi_{M_0}Q)(s)$ for $s\in\mathcal{A}_0$. Also, by
\eqref{Kcondgen}, we can
assume that
%
\begin{equation}
\label{jumpLBound}
K(n)\ge n^{-\alpha-1-\varepsilon}\qquad \forall n\geq M_0.
\end{equation}
By the ergodicity of $Q$, for every $s\in\mathcal{A}_0$ we have
%
\begin{equation}
\lim_{M\to\infty} \frac{1}{M} |\{0\leq j\le M-M_0\dvtx
\pi_{M_0}(\tilde\theta^j Y)=s\}| = f_s \qquad\mbox{for } Q\mbox{-a.e. } Y.\hspace*{-35pt}
\end{equation}
Consequently, there is a large $M$ and a finite set $\mathcal
{A}\subset
\widetilde E^M$ with
%
\begin{equation}
\label{typicalSentences}
(\pi_M Q)(\mathcal{A}) \geq1-\varepsilon
\end{equation}
such that
%
\begin{equation}
\label{frequencies}
\biggl|\frac{1}{M}|\{0\leq j\le M-M_0\dvtx \pi_{M_0}(\tilde\theta^j
z)=s\}|
-f_s\biggr| < \delta_1 \qquad\forall s\in\mathcal{A}_0, z\in\mathcal{A}.\hspace*{-35pt}
\end{equation}
Moreover, we can assume for all $z\in\mathcal{A}$ the following
relations, which are stated
in~\cite{BiGrdHo09}, equation (3.6), and are consequences of
ergodicity, too:
\begin{eqnarray}
&|\kappa(z)|  \in[M(m_Q-\varepsilon), M(m_Q+\varepsilon)],& \\
&\label{typicalConcatenation}
\hspace*{-98pt}\log Q\bigl(\kappa\bigl(Y^{(1)}, Y^{(2)}, \ldots, Y^{(M)}\bigr)=\kappa(z)\bigr)&
\nonumber
\\[-8pt]
\\[-8pt]
\nonumber
&\in\bigl[-M \bigl(m_Q H(\Psi_Q)+\varepsilon\bigr),-M \bigl(m_Q H(\Psi_Q)-\varepsilon\bigr)\bigr],&
\\
&\qquad\quad\log Q\bigl(\bigl(Y^{(1)}, Y^{(2)}, \ldots, Y^{(M)}\bigr)=z\bigr)
\in\bigl[-M \bigl(H(Q)+\varepsilon\bigr),-M\bigl (H(Q)-\varepsilon\bigr)\bigr],&
\label{typicalSentence}\\
&\displaystyle\sum_{i=1}^{|\kappa(z)|}\log\mu_0((\kappa(z))_i)-M m_Q \ev_{\Psi
_Q}[\log\mu_0(X_1)]
 \in[-M\varepsilon,M\varepsilon],&
\label{stringProbability}\\
&\displaystyle\sum_{i=1}^M \log K\bigl(\bigl|z^{(i)}\bigr|\bigr)-M \ev_Q[\log K(\tau_1)]
\in[-M\varepsilon, M\varepsilon].&
\label{wordProbability}
\end{eqnarray}
In the above relations, $|\kappa(z)|$ denotes the length of the string
$\kappa(z)$, $(\kappa(z))_i$
is the~$i$th letter of that string, $z^{(i)}$ is the~$i$th word of the
sentence $z$, $|z^{(i)}|$
is its length, while $H(Q), H(\Psi_Q)$ are the specific entropies of
the measures $Q, \Psi_Q$. In
the last relation, $\tau_1$ is distributed as the length of the first
word of an element of
$\widetilde{E}^\N$ drawn from $Q$. Finally, $M$ can be chosen such that
%
\begin{equation}
\label{MRestrictions}
M > \frac{8 M_0}{\delta_1},\qquad
\frac{1}{M} (\alpha+1+\varepsilon)\log[M(m_Q+\varepsilon
)+M_0]<\varepsilon.
\end{equation}


\subsection{Step 2: Good trajectories}
\label{B2}

For given $\omega\in E^{\N_0}$, we define a set $\mathcal
{T}_{\varepsilon, M}^\omega$ of trajectories
for the renewal sequence $T=(T_i)_{i\in\N_0}$ on which $R_N^\omega
\in
\mathcal{U}_{M, \delta_1}$.
In Step~3 we will control the probability that $T$ follows a trajectory in
$\mathcal{T}_{\varepsilon,M}^\omega$.

Let $\mathcal{B}:=\{\kappa(z): z\in\mathcal{A}\}$ be the set of
concatenations of the sentences of
$\mathcal{A}$. By~\eqref{typicalSentences} and \eqref{typicalConcatenation},
%
\begin{equation}
\label{stringNumber}
|\mathcal{B}|\ge(1-\varepsilon) \e^{M (m_Q H(\Psi_Q)-\varepsilon)}.
\end{equation}
Divide $\omega$ into consecutive pieces of length $\Lambda
:=[M(m_Q+\varepsilon)]+M_0$, mark with
1 those pieces that start with an element of $\mathcal{B}$, and mark
with 0 the remaining pieces,
that is, for $j\geq0$, let
%
\begin{equation}
\sigma_j := 1_{\{\theta^{j\Lambda}\omega\ \mathrm{starts\ with\ an\
element\
of\ } \mathcal{B}\}}.
\end{equation}
Let $\{j(r)\dvtx  r\geq1\}$ be the increasing sequence that picks out
the $j\geq1$ with
$\sigma_j=1$, and let $j(0)=0$. The increments $\{j(r+1)-j(r)\dvtx
r\geq0\}$ are i.i.d. geometric random variables with probability of
success $p_{\mathcal{B}}:=\mP(\omega\mbox{ starts
with an}\break \mbox{element of }\mathcal{B})$. It follows from \eqref
{stringProbability} and \eqref{stringNumber}
that
\begin{eqnarray}
p_{\mathcal{B}} &\geq&(1-\varepsilon)
\e^{M (m_Q H(\Psi_Q)+m_Q \ev_{\Psi_Q}[\log\mu
_0(X_1)]-2\varepsilon)}
\\
& =& (1-\varepsilon) \e^{-M m_Q H(\Psi_Q|\mu_0^{\otimes\N
_0})-2\varepsilon M}.
\label{ProbOfSuccess}
\end{eqnarray}
The equality in the second line follows from \cite{BiGrdHo09}, equation
(1.26). In particular, for
$\mP$-a.e. $\omega$ we have $\sigma_j=1$ for infinitely many $j$'s, and
so the sequence
$\{j(r)\dvtx  r\geq1\}$ is well defined.

Pick any $N>16 M/\delta_1$. The set $\mathcal{T}_{\varepsilon,
M}^\omega
$ consists of all $T$
that first jump to $j(1)\Lambda$ [i.e., $T_1=j(1)\Lambda$], next make
$M$ jumps that cut out
of $\theta^{j(1) \Lambda}\omega$ an element of $\mathcal{A}$ [which is
possible by the definitions
of $j(1)$ and $\mathcal{B}$], next jump to $j(2)\Lambda$ [i.e.,
$T_{M+2}=j(2)\Lambda$],
next again cut out an element of $\mathcal{A}$, and continue likewise
until they jump to
$j(\lceil N/(M+1)\rceil+1)\Lambda$ (no conditions are imposed
afterwards). The words between
two consecutive $j(r)\Lambda$'s we call a block. After the first jump
to $j(1)\Lambda$ and
up to the last jump to $j(\lceil N/(M+1)\rceil+1)\Lambda$, at least $N$
words are cut out,
because $T$ has created $\lceil N/(M+1)\rceil$ blocks each containing
exactly $M+1$ words.
We note that the first $M$ words are important and of typical length,
while the last word is
of an untypically large length and its sole purpose is for $T$ to move
to a good position
in $\omega$. Call $Y^{(1)},Y^{(2)},\ldots,Y^{(N)}$ the first $N$
words cut.

\begin{itlemma}
\label{claim}
$R_N^\omega\in\mathcal{U}_{M, \delta_1}$ for all $T$ in $\mathcal
{T}_{\varepsilon, M}^\omega$.
\end{itlemma}

\begin{pf}
By the definition of $R_N^\omega$, we need to show that every element
$s\in\mathcal{A}_0$
occurs in the finite sequence
%
\begin{equation}
\label{finiteSequence}
\bigl(\pi_{M_0} \tilde\theta^j \bigl(Y^{(1)},Y^{(2)},\ldots,Y^{(N)}\bigr)^{\mathrm
{per}}\bigr)_{0\leq j\leq N-1}
\end{equation}
the right number of times, that is, a number of times that falls in the
interval $((f_s-2\delta_1)N,
(f_s+2\delta_1) N)$.\vadjust{\goodbreak}

For the lower bound, note that the sequence $(Y^{(1)},Y^{(2)},\ldots
,Y^{(N)})$ contains at least
the words of the first $[N/(M+1)]-1$ blocks out of the $\lceil
N/(M+1)\rceil$ blocks that $T$
created, because the last word of these blocks has index at most
$i^*=1+(N/(M+1)-1)(M+1)=N-M<N$.
Each such block offers at least $M (f_s-\delta_1)$ occurrences of the
word $s$, because of
\eqref{frequencies} and $i^*<N-M_0$. Thus, we have at least
%
\begin{eqnarray}
\qquad M(f_s-\delta_1)\biggl(\frac{N}{M+1}-2\biggr)
&=& N (f_s-\delta_1)-\frac{f_s-\delta_1}{M+1} N-2 M( f_s-\delta_1)
\nonumber
\\[-8pt]
\\[-8pt]
\nonumber
&>& N (f_s-2\delta_1)
\end{eqnarray}
occurrences of $s$ in the sequence in \eqref{finiteSequence}, where the
last inequality holds
because $N>16M/\delta_1$ and $M>8M_0/\delta_1$ by \eqref{MRestrictions}.

For the upper bound, note that, because of \eqref{frequencies}, the
occurrences of $s$ in the
sequence in \eqref{finiteSequence} are at most
%
\begin{eqnarray}
\qquad&& 1+\bigl((f_s+\delta_1)M+M_0\bigr)\biggl(\frac{N}{M+1}+1\biggr)+M_0
\nonumber
\\[-8pt]
\\[-8pt]
\nonumber
&&\qquad \leq N(f_s+\delta_1)+\frac{N M_0}{M+1}+M(f_s+\delta_1)+2M_0+1
< N (f_s+2\delta_1),
\end{eqnarray}
where the last inequality again uses $N>16M/\delta_1$ and
$M>8M_0/\delta_1$.
\end{pf}


\subsection{Step 3: Probability of good trajectories}
\label{B3}

For the quenched probability $\pr(T\in\mathcal{T}_{\varepsilon,
M}^\omega)$, we have the lower bound
%
\begin{eqnarray}
\label{ProbabilityLB}
\qquad\pr(T\in\mathcal{T}_{\varepsilon, M}^\omega) &
\geq& K(j(1)\Lambda)\nonumber\\
&&{} \times\bigl(\e^{M(H(Q)-m_Q H(\Psi_Q)-2\varepsilon)}
\e^{M (\ev_Q[\log K(\tau_1)]-\varepsilon)}\bigr)^{\lceil N/(M+1)\rceil
}
\\
&&{} \times\prod_{r=1}^{\lceil N/(M+1)\rceil}
\inf_{|\eta-M m_Q|<M \varepsilon}
K\bigl([j(r+1)-j(r)]\Lambda-\eta\bigr).\nonumber
\end{eqnarray}
The last product is a lower bound for the probability of the large
jumps that land at the
points $j(r+1)\Lambda$, $1\leq r\leq\lceil N/(M+N)\rceil$. The power
preceding this product
corresponds to the jumps inside each of the $\lceil N/(M+1)\rceil$
blocks, and uses that,
by (\ref{typicalConcatenation}) and (\ref{typicalSentence}), for each
element of $\mathcal{B}$
there are at least $\e^{M(H(Q)-m_Q H(\Psi_Q)-2\varepsilon)}$ different
words of $\mathcal{A}$
having this element as concatenation, and that, by \eqref
{wordProbability}, the probability
for $M$ jumps to cut out a given word in $\mathcal{A}$ is at least $\e
^{M (\mE_Q[\log K(\tau_1)]
-\varepsilon)}$. It therefore follows that
%
\begin{eqnarray}
\label{precise}
&&\liminf_{N\to\infty} \frac{1}{N} \log\pr(\tau\in\mathcal
{T}_{\varepsilon, M}^\omega)\nonumber\\
&&\qquad\geq H(Q)-m_Q H(\Psi_Q)+\ev_Q[\log K(\tau_1)]-3\varepsilon
\\
&&\qquad\quad{}+ \frac{1}{M} \mE\Bigl(\log\Bigl[\inf_{|\eta-M m_Q|<M \varepsilon_1}
K\bigl([j(2)-j(1)]\Lambda-\eta\bigr)\Bigr]\Bigr).\nonumber
\end{eqnarray}
To be more precise, \eqref{ProbabilityLB} gives \eqref{precise} with
the right-hand side multiplied
by $M/(M+1)$, but since the factors in \eqref{ProbabilityLB} are
probabilities, replacing $M/(M+1)$
by 1 still gives us a lower bound. Now, because of \eqref{jumpLBound}
and $\Lambda-\eta\ge M_0$,
the last expectation is bounded from below by
%
\begin{eqnarray}
\label{Elogs}
&&\mE\bigl[\log\bigl(\bigl([j(2)-j(1)]\Lambda\bigr)^{-\alpha-1-\varepsilon}\bigr)\bigr]\nonumber\\
&&\qquad= -(\alpha+1+\varepsilon) \mE\bigl[\log\bigl([j(2)-j(1)]\Lambda\bigr)\bigr]
\\
&&\qquad\geq-(\alpha+1+\varepsilon)\bigl(\log\Lambda+\log\mE[j(2)-j(1)]\bigr),\nonumber
\end{eqnarray}
where we use the concavity of $\log$. Since $\mE
[j(2)-j(1)]=1/p_\mathcal
{B}$, by combining
(\ref{precise}) and (\ref{Elogs}) with the lower bound on $p_\mathcal{B}$ in
\eqref{ProbOfSuccess},
we get that
%
\begin{eqnarray}
&&\liminf_{N\to\infty} N^{-1}\log\pr(T\in\mathcal{T}_{\varepsilon,
M}^\omega)\nonumber\\
&&\qquad \geq H(Q)-m_Q H(\Psi_Q)+\ev_Q[\log K(\tau_1)]-3\varepsilon\nonumber\\
&&\qquad\quad{} + \frac{1}{M} \bigl[-(\alpha+1+\varepsilon)\nonumber\\
&&\hspace*{39pt}\qquad{}\times \bigl(\log\Lambda-\log
(1-\varepsilon)
+ M m_Q H(\Psi_Q|\mu_0^{\otimes\N_0}) + 2\varepsilon M \bigr)\bigr]
\\
&&\qquad = H(Q) - m_Q H(\Psi_Q) + \ev_Q[\log K(\tau_1)]
- m_Q H(\Psi_Q|\mu_0^{\otimes\N_0})\nonumber\\
&&\qquad\quad{} -\alpha m_Q H(\Psi_Q|\mu
_0^{\otimes
\N_0}) -3\varepsilon-\varepsilon m_Q H(\Psi_Q|\mu_0^{\otimes\N_0})
\nonumber\\
&&\qquad\quad{}-\frac{1}{M}(\alpha+1+\varepsilon)\bigl(\log\Lambda-\log(1-\varepsilon)\bigr)-2 (\alpha+1+\varepsilon)\varepsilon.\nonumber
\end{eqnarray}
The fourth line equals $-I^\mathrm{que}(Q)$ because of \cite
{BiGrdHo09}, equations (1.16), (1.30)
and (1.32), where in using (1.16) we note that what we call in this
paper $\alpha$ is called $\alpha-1$ in \cite{BiGrdHo09}.
The fifth line is at least $-\varepsilon C_Q$ for some positive
constant $C_Q$ that
depends on $Q$, because of \eqref{MRestrictions}. Thus, we end up with
%
\begin{equation}
\label{probabilityLowerBound}
\liminf_{N\to\infty} N^{-1}\log\pr(T\in\mathcal{T}_{\varepsilon,
M}^\omega)
\geq-I^\mathrm{que}(Q)-\varepsilon C_Q.
\end{equation}


\subsection{Step 4: Lower bound}
\label{B4}

For $T\in\mathcal{T}_{\varepsilon, M}^\omega$, we have
\begin{equation} \label{LengthUpperBound}
N m_{R_N^\omega} \leq j\biggl(\biggl\lceil\frac{N}{M+1}\biggr\rceil+1\biggr)
\bigl(M(m_Q+\varepsilon)+M_0\bigr),
\end{equation}
and note that
\begin{equation}
\lim_{N\to\infty} \frac{1}{N}j\biggl(\biggl\lceil\frac{N}{M+1}\biggr\rceil+1\biggr)
= \frac{1}{M+1}\frac{1}{p_\mathcal{B}}. \label{LawOfLN}
\end{equation}
Hence
%
\begin{eqnarray}
\label{SLBound}
&&\ev\bigl(\e^{N m_{R_N^\omega} \log z +
N\beta\Phi(R_N^\omega)}\bigr)\nonumber\\
&&\qquad\geq\ev\bigl(\e^{N m_{R_N^\omega}\log z + N\beta\Phi(R_N^\omega)}
1_{\{T\in\mathcal{T}_{\varepsilon, M}^\omega\}}\bigr)
\\
&&\qquad\geq\e^{N \beta(\Phi(Q)-\varepsilon)}
z^{j(\lceil N/(M+1)\rceil+1)(M(m_Q+\varepsilon)+M_0)}
\pr(\tau\in\mathcal{T}_{\varepsilon, M}^\omega).\nonumber
\end{eqnarray}
Combining \eqref{Sdefinition}, \eqref{probabilityLowerBound}, \eqref
{LawOfLN} and \eqref{SLBound}, we get
%
\begin{eqnarray}
S^\mathrm{que}(\beta;z) &\geq&\beta\Phi(Q)-\beta\varepsilon
+ \frac{M(m_Q+\varepsilon)+M_0}{(M+1)p_\mathcal{B}} \log z
\nonumber
\\[-8pt]
\\[-8pt]
\nonumber
&&{}-I^\mathrm{que}(Q)-C_Q \varepsilon.
\end{eqnarray}
Now let $z \uparrow1$ and $\varepsilon\downarrow0$, to get
$S^\mathrm
{que}(\beta;1-)
\geq\beta\Phi(Q)-I^\mathrm{que}(Q)>A$. Since $A<S^\mathrm
{que}(\beta
;1)$ was arbitrary,
it follows that $S^\mathrm{que}(\beta;1-) \ge S^\mathrm{que}(\beta;1)$.


\end{appendix}

\section*{Acknowledgment}
The research in this paper was carried out while the first author was
a postdoc at EURANDOM.

%

\printaddresses

\end{document}